\newcommand{\xBc}{\langle}
\newcommand{\xBe}{\rangle}

\newcommand{\xbD}{\Delta}

\newcommand{\xbG}{\Gamma}

\newcommand{\xbP}{\Pi}

\newcommand{\xbS}{\Sigma}

\newcommand{\xba}{\alpha}
\newcommand{\xbb}{\beta}

\newcommand{\xbe}{\in}
\newcommand{\xbf}{\phi}
\newcommand{\xbg}{\gamma}

\newcommand{\xbm}{\mu}

\newcommand{\xbo}{\omega}

\newcommand{\xbq}{\psi}
\newcommand{\xbr}{\rho}
\newcommand{\xbs}{\sigma}
\newcommand{\xbt}{\tau}

\newcommand{\xCK}{\times}

\newcommand{\xCN}{\neg}
\newcommand{\xCO}{ }
\newcommand{\xCQ}{\emptyset}

\newcommand{\xCd}{\approx}

\newcommand{\xCf}{\hspace{0.1em}}

\newcommand{\xCq}{\sim}

\newcommand{\xcA}{\forall}

\newcommand{\xcC}{\not\subseteq}

\newcommand{\xcE}{\exists}

\newcommand{\xcH}{\not\Rightarrow}
\newcommand{\xcI}{\not\Leftarrow}

\newcommand{\xcL}{\not\vdash}

\newcommand{\xcN}{\hspace{0.2em}\not\sim\hspace{-0.9em}\mid\hspace{0.8em}}

\newcommand{\xcT}{\bot}

\newcommand{\xcb}{\subset}
\newcommand{\xcc}{\subseteq}
\newcommand{\xcd}{\supseteq}
\newcommand{\xce}{\not\in}

\newcommand{\xcg}{\geq}
\newcommand{\xch}{\Rightarrow}
\newcommand{\xci}{\Leftarrow}
\newcommand{\xcj}{\Leftrightarrow}
\newcommand{\xck}{\leq}
\newcommand{\xcl}{\vdash}
\newcommand{\xcm}{\models}
\newcommand{\xcn}{\hspace{0.2em}\sim\hspace{-0.9em}\mid\hspace{0.58em}}

\newcommand{\xco}{\vee}
\newcommand{\xcp}{\rightarrow}

\newcommand{\xcr}{\leftrightarrow}
\newcommand{\xcs}{\cap}
\newcommand{\xcu}{\wedge}
\newcommand{\xcv}{\cup}

\newcommand{\xcz}{\Box}

\newcommand{\xDH}{\item }

\newcommand{\xDO}{\circ}

\newcommand{\xdD}{\mbox{\boldmath$D$}}

\newcommand{\xdc}{{\cal C}}
\newcommand{\xdd}{{\cal D}}

\newcommand{\xdf}{{\cal F}}

\newcommand{\xdi}{{\cal I}}

\newcommand{\xdl}{{\cal L}}
\newcommand{\xdm}{{\cal M}}

\newcommand{\xdp}{{\cal P}}

\newcommand{\xdy}{{\cal Y}}

\newcommand{\xEH}{ & }
\newcommand{\xEI}{\begin{itemize}}
\newcommand{\xEJ}{\end{itemize}}
\newcommand{\xEP}{ \\ }

\newcommand{\xEc}{\not<}
\newcommand{\xEd}{\neq}

\newcommand{\xEh}{\begin{enumerate}}

\newcommand{\xEj}{\end{enumerate}}

\newcommand{\xeA}{\nabla}
\newcommand{\xeB}{\not\prec}

\newcommand{\xeb}{\prec}
\newcommand{\xec}{\preceq}

\newcommand{\xex}{\upharpoonright}

\newcommand{\xFO}{\parallel}

\newcommand{\xfA}{\mid}

\newcommand{\Xl}{\ldots}

\newcommand{\ol}{\overline}

\newcommand{\xssc}{\scriptsize}

\newcommand{\xssB}{\scriptsize}

\newcommand{\bl}{\begin{lemma} \rm}
\newcommand{\el}{\end{lemma}}
\newcommand{\br}{\begin{remark} \rm}
\newcommand{\er}{\end{remark}}
\newcommand{\be}{\begin{example} \rm}
\newcommand{\ee}{\end{example}}
\newcommand{\bco}{\begin{corollary} \rm}
\newcommand{\eco}{\end{corollary}}
\newcommand{\bc}{\begin{claim} \rm}
\newcommand{\ec}{\end{claim}}
\newcommand{\bfa}{\begin{fact} \rm}
\newcommand{\efa}{\end{fact}}
\newcommand{\bp}{\begin{proposition} \rm}
\newcommand{\ep}{\end{proposition}}
\newcommand{\bd}{\begin{definition} \rm}
\newcommand{\ed}{\end{definition}}
\newcommand{\bcs}{\begin{construction} \rm}
\newcommand{\ecs}{\end{construction}}
\newcommand{\bcd}{\begin{condition} \rm}
\newcommand{\ecd}{\end{condition}}
\newcommand{\bt}{\begin{theorem} \rm}
\newcommand{\et}{\end{theorem}}
\newcommand{\bn}{\begin{notation} \rm}
\newcommand{\en}{\end{notation}}
\newcommand{\bfi}{\begin{bild} \rm}
\newcommand{\efi}{\end{bild}}
\newcommand{\bsta}{\begin{statement} \rm}
\newcommand{\esta}{\end{statement}}
\newcommand{\bcom}{\begin{comment} \rm}
\newcommand{\ecom}{\end{comment}}
\newcommand{\bdia}{\begin{diagram} \rm}
\newcommand{\edia}{\end{diagram}}

\newcommand{\bfc}{\begin{figure}[htb] \begin{center}}
\newcommand{\efc}{\end{center} \end{figure}}

\documentclass{article}

\usepackage{amssymb,latexsym,epic,eepic,rotating}

\title{
Independence and abstract multiplication
\thanks{
paper 367
}
}

\author{Dov M Gabbay
\thanks{
Dov.Gabbay@kcl.ac.uk, www.dcs.kcl.ac.uk/staff/dg
} \\
King's College, London
\thanks{
Department of Computer Science, King's College London, Strand,
London WC2R 2LS, UK
} \\
and \\
Bar-Ilan University, Israel
\thanks{
Department of Computer Science,
Bar-Ilan University,
52900 Ramat-Gan, Israel
} \\
and \\
University of Luxembourg
\thanks{
Computer Science and Communications,
Faculty of Sciences,
6, rue Coudenhove-Kalergi,
L-1359 Luxembourg
} \\ \\
Karl Schlechta
\thanks{
ks@cmi.univ-mrs.fr, karl.schlechta@web.de, http://www.cmi.univ-mrs.fr/ $\sim$ ks
} \\
Laboratoire d'Informatique Fondamentale de Marseille
\thanks{
UMR 6166, CNRS and Universit\'{e} de Provence,
Address: CMI, 39, rue Joliot-Curie, F-13453 Marseille Cedex 13, France
}
}

\begin{document}

\newtheorem{lemma}{Lemma}[section]
\newtheorem{theorem}[lemma]{Theorem}
\newtheorem{proposition}[lemma]{Proposition}
\newtheorem{corollary}[lemma]{Corollary}
\newtheorem{claim}[lemma]{Claim}
\newtheorem{fact}[lemma]{Fact}
\newtheorem{remark}[lemma]{Remark}
\newtheorem{definition}{Definition}[section]
\newtheorem{construction}{Construction}[section]
\newtheorem{condition}{Condition}[section]
\newtheorem{example}{Example}[section]
\newtheorem{notation}{Notation}[section]
\newtheorem{bild}{Figure}[section]
\newtheorem{comment}{Comment}[section]
\newtheorem{statement}{Statement}[section]
\newtheorem{diagram}{Diagram}[section]

\renewcommand{\labelenumi}
  {(\arabic{enumi})}
\renewcommand{\labelenumii}
  {(\arabic{enumi}.\arabic{enumii})}
\renewcommand{\labelenumiii}
  {(\arabic{enumi}.\arabic{enumii}.\arabic{enumiii})}
\renewcommand{\labelenumiv}
  {(\arabic{enumi}.\arabic{enumii}.\arabic{enumiii}.\arabic{enumiv})}

\maketitle

\setcounter{secnumdepth}{3}
\setcounter{tocdepth}{3}

\begin{abstract}

We investigate the notion of independence, which is at the basis of many,
seemingly unrelated, properties of logics, like the Rational Monotony rule
of nonmonotonic logics, but also of interpolation theorems of monotonic and
nonmonotonic logic. We show a strong connection between independence and
certain rules about multiplication of abstract size in the field of
nonmonotonic logic. We think that this notion of independence, with its
ramifications, is extremely important, and has not been sufficiently
investigated.

\end{abstract}

\tableofcontents

%
%
%


\section{
Abstract definition of independence
}

\label{Section Mul-Def}

$ \xCO $

\vspace{10mm}

\begin{diagram}

\label{Diagram Mul-Commut}
\index{Diagram Mul-Commut}

\centering
\setlength{\unitlength}{1mm}
{\renewcommand{\dashlinestretch}{30}
\begin{picture}(150,100)(0,0)

\path(35,90)(78,90)
\path(75.2,91)(78,90)(75.2,89)
\path(66,60)(106,60)
\path(103.2,61)(106,60)(103.2,59)
\path(35,30)(78,30)
\path(75.2,31)(78,30)(75.2,29)

\path(35,35)(55,55)
\path(52.3,53.7)(55,55)(53.7,52.3)
\path(88,35)(108,55)
\path(105.3,53.7)(108,55)(106.7,52.3)
\path(90.7,36.3)(88,35)(89.3,37.7)

\path(35,85)(55,65)
\path(53.7,67.7)(55,65)(52.3,66.3)
\path(88,85)(108,65)
\path(106.7,67.7)(108,65)(105.3,66.3)
\path(89.3,82.3)(88,85)(90.7,83.7)

\put(29,89){{\xssc $\xbS_1$}}
\put(29,29){{\xssc $\xbS_2$}}
\put(55,59){{\xssc $\xbS_1 \xDO \xbS_2$}}

\put(81,89){{\xssc $f(\xbS_1$)}}
\put(81,29){{\xssc $f(\xbS_2$)}}
\put(108,59){{\xssc $f(\xbS_1 \xDO \xbS_2)=f(\xbS_1) \xDO' f(\xbS_2$)}}

\put(35,10){{\xssc Note that $\xDO$ and $\xDO'$ might be different}}
\put(50,1){{\xssc Independence}}

\end{picture}
}

\end{diagram}

\vspace{4mm}

$ \xCO $

The right notion of independence in our context seems to be:

We have compositions $ \xDO $ and $ \xDO',$ and operation $f.$ We can
calculate
$f(\xbS_{1} \xDO \xbS_{2})$ from $f(\xbS_{1})$ and $f(\xbS_{2}),$ but
also conversely, given
$f(\xbS_{1} \xDO \xbS_{2})$ we can calculate $f(\xbS_{1})$ and $f(
\xbS_{2}).$
Of course, in other contexts, other notions of independence
might be adequate.
More precisely:

\bd

$\hspace{0.01em}$


\label{Definition Mul-Ind}

Let $f: \xdd \xcp \xdc $ be any function from domain $ \xdd $ to co-domain
$ \xdc.$
Let $ \xDO $ be a ``composition function'' $ \xDO: \xdd \xCK \xdd \xcp \xdd
,$ likewise
for $ \xDO': \xdc \xCK \xdc \xcp \xdc.$

We say that $ \xBc f, \xDO, \xDO'  \xBe $ are independent iff for any $
\xbS_{i}
\xbe \xdd $

(1) $f(\xbS_{1} \xDO \xbS_{2})=f(\xbS_{1}) \xDO' f(\xbS_{2}),$

(2) we can recover $f(\xbS_{i})$ from $f(\xbS_{1} \xDO \xbS_{2}),$
provided we know how
$ \xbS_{1} \xDO \xbS_{2}$ splits into the $ \xbS_{i},$ without using $f$
again.
\subsection{
Discussion
}

\ed

 \xEh

 \xDH Ranked structures satisfy it:

Let $ \xDO = \xDO' = \xcv.$ Let $f$ be the minimal model operator $ \xbm
$ of preferential logic.
Let $X,Y \xcc X \xcv Y$ have (at least) medium size (see below).
Then $ \xbm (X \xcv Y)= \xbm (X) \xcv \xbm (Y),$ and $ \xbm (X)= \xbm (X
\xcv Y) \xcs X,$ $ \xbm (Y)= \xbm (X \xcv Y) \xcs Y.$

 \xDH Consistent classical formulas and their interpretation satisfy it:

Let $ \xDO $ be conjunction in the composed language,
$ \xDO' $ be model set intersection, $f(\xbf)=M(\xbf).$
Let $ \xbf,$ $ \xbq $ be classical formulas, defined on disjoint language
fragments
$ \xdl,$ $ \xdl' $ of some language $ \xdl''.$ Then $f(\xbf \xcu \xbq
)=M(\xbf) \xcs M(\xbq),$ and
$M(\xbf)$ is the projection of $M(\xbf) \xcs M(\xbq)$ onto the
(models of) language $ \xdl,$
likewise for $M(\xbq).$
This is due to the way validity is defined, using only variables which
occur in the formula.

As a consequence, monotonic logic has semantical interpolation -
see  \cite{GS09c}, and below,
Section \ref{Section Mul-Mon-Int} (page \pageref{Section Mul-Mon-Int}). The
definition of being insensitive
is justified by
this modularity.

 \xDH It does not hold for inconsistent classical formulas:
We cannot recover $M(a \xcu \xCN a)$ and $M(b)$ from $M(a \xcu \xCN a \xcu
b),$ as we do not
know where the inconsistency came from.
The basic reason is trivial: One empty factor suffices to make the whole
product empty, and we do not know which factor was the culprit.
See Section \ref{Section Mul-Relev} (page \pageref{Section Mul-Relev})  for the
discussion of a remedy.

 \xDH Preferential logic satisfies it under certain conditions:

If $ \xbm (X \xCK Y)= \xbm (X) \xCK \xbm (Y)$ holds for model products and
$ \xcn,$ then it
holds by definition. An important consequence is that such a logic
has interpolation of the form $ \xcn \xDO \xcn,$
see Section \ref{Section Mul-Int} (page \pageref{Section Mul-Int}).

 \xDH Modular revision a la Parikh is based on a similar idea.

 \xEj
\subsection{
Independence and multiplication of abstract size
}

We are mainly interested in nonmonotonic logic. In this domain,
independence is strongly connected to multiplication of
abstract size, and much of the present paper treats this connection
and its repercussions.

We have at least two scenarios for multiplication, one is decribed in
Diagram \ref{Diagram Mul-Add} (page \pageref{Diagram Mul-Add}), the second in
Diagram \ref{Diagram Mul-Prod} (page \pageref{Diagram Mul-Prod}).
In the first scenario, we have nested sets, in the second, we have set
products.
In the first scenario, we consider subsets which behave as the big set
does,
in the second scenario we consider subspaces, and decompose the behaviour
of the big space into behaviour of the subspaces.
In both cases, this results naturally in multiplication of abstract sizes.
When we look at the corresponding relation properties, they are quite
different (rankedness vs. some kind of modularity). But this is perhaps
to be expected, as the two scenarios are quite different.

We do not know whether there are still other,
interesting, scenarios to consider in our framework.
\section{
Introduction to abstract size, additive rules
}

\label{Section Mul-Intro}

To put our work more into perspective, we first repeat in this section
material from  \cite{GS08c}. This gives the main definitions and
rules
for non-monotonic logics, see
Table \ref{Table Base2-Rules-Def-Conn-1} (page \pageref{Table
Base2-Rules-Def-Conn-1})  and
Table \ref{Table Base2-Rules-Def-Conn-2} (page \pageref{Table
Base2-Rules-Def-Conn-2}),
``Logical rules, definitions and connections''.
We then give the main additive rules for
manipulation of abstract size from  \cite{GS09a}, see
Table \ref{Table Base2-Size-Rules-1} (page \pageref{Table Base2-Size-Rules-1}) 
and
Table \ref{Table Base2-Size-Rules-2} (page \pageref{Table Base2-Size-Rules-2}),
``Rules on size''.

$ \xCO $

\label{Definition Log-Cond-NoRef-Size}
\index{Definition Logical conditions}

Explanation of
Table \ref{Table Base2-Rules-Def-Conn-1} (page \pageref{Table
Base2-Rules-Def-Conn-1}),
``Logical rules, definitions and connections Part $I'' $
and
Table \ref{Table Base2-Rules-Def-Conn-2} (page \pageref{Table
Base2-Rules-Def-Conn-2}),
''Logical rules, definitions and connections Part II":

The tables are split in two, as they would not fit onto a page otherwise.
The difference between the first two columns is that the first column
treats the formula version of the rule, the second the more general
theory (i.e., set of formulas) version.

The first column ``Corr.''
is to be understood as follows:

Let a logic $ \xcn $ satisfy $ \xCf (LLE)$ and $ \xCf (CCL),$ and define a
function $f: \xdD_{ \xdl } \xcp \xdD_{ \xdl }$
by $f(M(T)):=M(\ol{ \ol{T} }).$ Then $f$ is well defined, satisfies $(
\xbm dp),$ and $ \ol{ \ol{T} }=Th(f(M(T))).$

If $ \xcn $ satisfies a rule in the left hand side,
then - provided the additional properties noted in the middle for $ \xch $
hold, too -
$f$ will satisfy the property in the right hand side.

Conversely, if $f: \xdy \xcp \xdp (M_{ \xdl })$ is a function, with $
\xdD_{ \xdl } \xcc \xdy,$ and we define a logic
$ \xcn $ by $ \ol{ \ol{T} }:=Th(f(M(T))),$ then $ \xcn $ satisfies $ \xCf
(LLE)$ and $ \xCf (CCL).$
If $f$ satisfies $(\xbm dp),$ then $f(M(T))=M(\ol{ \ol{T} }).$

If $f$ satisfies a property in the right hand side,
then - provided the additional properties noted in the middle for $ \xci $
hold, too -
$ \xcn $ will satisfy the property in the left hand side.

We use the following abbreviations for those supplementary conditions in
the
``Correspondence'' columns:
``$T= \xbf $'' means that, if one of the theories
(the one named the same way in Definition \ref{Definition Log-Cond-NoRef-Size}
(page \pageref{Definition Log-Cond-NoRef-Size}))
is equivalent to a formula, we do not need $(\xbm dp).$
$-(\xbm dp)$ stands for
``without $(\xbm dp)$''.

$A=B \xFO C$ will abbreviate $A=B,$ or $A=C,$ or $A=B \xcv C.$

\begin{table}[h]

\index{$(Opt)$}
\index{$(\xbm \xcc)$}
\index{$(REF)$}
\index{Reflexivity}
\index{$(LLE)$}
\index{Left Logical Equivalence}
\index{$(RW)$}
\index{Right Weakening}
\index{$(iM)$}
\index{$(wOR)$}
\index{$(\xbm wOR)$}
\index{$(eM\xdi)$}
\index{$(disjOR)$}
\index{$(\xbm disjOR)$}
\index{$(I\xcv disj)$}
\index{$(CP)$}
\index{$(\xbm \xCQ)$}
\index{$(I_1)$}
\index{Consistency Preservation}
\index{$(\xbm \xCQ fin)$}
\index{$(AND_1)$}
\index{$(I_2)$}
\index{$(AND_n)$}
\index{$(I_n)$}
\index{$(AND)$}
\index{$(I_\xbo)$}
\index{$(CCL)$}
\index{Classical Closure}
\index{$(iM)$}
\index{$(OR)$}
\index{$(\xbm OR)$}
\index{$(PR)$}
\index{$(\xbm PR)$}
\index{$(\xbm PR')$}
\index{$(CUT)$}
\index{$ (\xbm CUT) $}

\tabcolsep=0.5pt
\caption{Logical rules, definitions and connections Part I}

\label{Table Base2-Rules-Def-Conn-1}
\begin{center}
\begin{turn}{90}
{\tiny
\begin{tabular}{|c|c|c|c|c|c|}

\hline

\multicolumn{6}{|c|}{\bf Logical rules, definitions and connections Part I}\\
\hline

\multicolumn{2}{|c|}{Logical rule}
\xEH
Corr.
\xEH
Model set
\xEH
Corr.
\xEH
Size Rules
\xEP

\hline

\multicolumn{6}{|c|}{Basics}
\xEP

\hline

$(SC)$ Supraclassicality
\index{$(SC)$}
\index{Supraclassicality}
\xEH
$(SC)$
\xEH
$\xch$
\xEH
$(\xbm \xcc)$
\xEH
trivial
\xEH
$(Opt)$
\xEP

\cline{3-3}

$ \xba \xcl \xbb $ $ \xch $ $ \xba \xcn \xbb $
\xEH
$ \ol{T} \xcc \ol{ \ol{T} }$
\xEH
$\xci$
\xEH
$f(X) \xcc X$
\xEH
\xEH
\xEP

\cline{1-1}

$(REF)$ Reflexivity
\xEH
\xEH
\xEH
\xEH
\xEH
\xEP

$ T \xcv \{\xba\} \xcn \xba $
\xEH
\xEH
\xEH
\xEH
\xEH
\xEP

\hline

$(LLE)$
\xEH
$(LLE)$
\xEH
\xEH
\xEH
\xEH
\xEP

Left Logical Equivalence
\xEH
\xEH
\xEH
\xEH
\xEH
\xEP

$ \xcl \xba \xcr \xba',  \xba \xcn \xbb   \xch $
\xEH
$ \ol{T}= \ol{T' }  \xch   \ol{\ol{T}} = \ol{\ol{T'}}$
\xEH
\xEH
\xEH
\xEH
\xEP

$ \xba' \xcn \xbb $
\xEH
\xEH
\xEH
\xEH
\xEH
\xEP

\hline

$(RW)$ Right Weakening
\xEH
$(RW)$
\xEH
\xEH
\xEH
trivial
\xEH
$(iM)$
\xEP

$ \xba \xcn \xbb,  \xcl \xbb \xcp \xbb'   \xch $
\xEH
$ T \xcn \xbb,  \xcl \xbb \xcp \xbb'   \xch $
\xEH
\xEH
\xEH
\xEH
\xEP

$ \xba \xcn \xbb' $
\xEH
$T \xcn \xbb' $
\xEH
\xEH
\xEH
\xEH
\xEP

\hline

$(wOR)$
\xEH
$(wOR)$
\xEH
$\xch$
\xEH
$(\xbm wOR)$
\xEH
$\xcj$
\xEH
$(eM\xdi)$
\xEP

\cline{3-3}

$ \xba \xcn \xbb,$ $ \xba' \xcl \xbb $ $ \xch $
\xEH
$ \ol{ \ol{T} } \xcs \ol{T' }$ $ \xcc $ $ \ol{ \ol{T \xco T' } }$
\xEH
$\xci$
\xEH
$f(X \xcv Y) \xcc f(X) \xcv Y$
\xEH
\xEH
\xEP

$ \xba \xco \xba' \xcn \xbb $
\xEH
\xEH
\xEH
\xEH
\xEH
\xEP

\hline

$(disjOR)$
\xEH
$(disjOR)$
\xEH
$\xch$
\xEH
$(\xbm disjOR)$
\xEH
$\xcj$
\xEH
$(I\xcv disj)$
\xEP

\cline{3-3}

$ \xba \xcl \xCN \xba',$ $ \xba \xcn \xbb,$
\xEH
$\xCN Con(T \xcv T') \xch$
\xEH
$\xci$
\xEH
$X \xcs Y= \xCQ $ $ \xch $
\xEH
\xEH
\xEP

$ \xba' \xcn \xbb $ $ \xch $ $ \xba \xco \xba' \xcn \xbb $
\xEH
$ \ol{ \ol{T} } \xcs \ol{ \ol{T' } } \xcc \ol{ \ol{T \xco T' } }$
\xEH
\xEH
$f(X \xcv Y) \xcc f(X) \xcv f(Y)$
\xEH
\xEH
\xEP

\hline

$(CP)$
\xEH
$(CP)$
\xEH
$\xch$
\xEH
$(\xbm \xCQ)$
\xEH
trivial
\xEH
$(I_1)$
\xEP

\cline{3-3}

Consistency Preservation
\xEH
\xEH
$\xci$
\xEH
\xEH
\xEH
\xEP

$ \xba \xcn \xcT $ $ \xch $ $ \xba \xcl \xcT $
\xEH
$T \xcn \xcT $ $ \xch $ $T \xcl \xcT $
\xEH
\xEH
$f(X)= \xCQ $ $ \xch $ $X= \xCQ $
\xEH
\xEH
\xEP

\hline

\xEH
\xEH
\xEH
$(\xbm \xCQ fin)$
\xEH
\xEH
$(I_1)$
\xEP

\xEH
\xEH
\xEH
$X \xEd \xCQ $ $ \xch $ $f(X) \xEd \xCQ $
\xEH
\xEH
\xEP

\xEH
\xEH
\xEH
for finite $X$
\xEH
\xEH
\xEP

\hline

\xEH
$(AND_1)$
\xEH
\xEH
\xEH
\xEH
$(I_2)$
\xEP

\xEH
$\xba\xcn\xbb \xch \xba\xcN\xCN\xbb$
\xEH
\xEH
\xEH
\xEH
\xEP

\hline

\xEH
$(AND_n)$
\xEH
\xEH
\xEH
\xEH
$(I_n)$
\xEP

\xEH
$\xba\xcn\xbb_1, \ldots, \xba\xcn\xbb_{n-1} \xch $
\xEH
\xEH
\xEH
\xEH
\xEP

\xEH
$\xba\xcN(\xCN\xbb_1 \xco \ldots \xco \xCN\xbb_{n-1})$
\xEH
\xEH
\xEH
\xEH
\xEP

\hline

$(AND)$
\xEH
$(AND)$
\xEH
\xEH
\xEH
trivial
\xEH
$(I_\xbo)$
\xEP

$ \xba \xcn \xbb,  \xba \xcn \xbb'   \xch $
\xEH
$ T \xcn \xbb, T \xcn \xbb'   \xch $
\xEH
\xEH
\xEH
\xEH
\xEP

$ \xba \xcn \xbb \xcu \xbb' $
\xEH
$ T \xcn \xbb \xcu \xbb' $
\xEH
\xEH
\xEH
\xEH
\xEP

\hline

$(CCL)$ Classical Closure
\xEH
$(CCL)$
\xEH
\xEH
\xEH
trivial
\xEH
$(iM)+(I_\xbo)$
\xEP

\xEH
$ \ol{ \ol{T} }$ classically closed
\xEH
\xEH
\xEH
\xEH
\xEP

\hline

$(OR)$
\xEH
$(OR)$
\xEH
$\xch$
\xEH
$(\xbm OR)$
\xEH
$\xcj$
\xEH
$(eM\xdi)+(I_\xbo)$
\xEP

\cline{3-3}

$ \xba \xcn \xbb,  \xba' \xcn \xbb   \xch $
\xEH
$ \ol{\ol{T}} \xcs \ol{\ol{T'}} \xcc \ol{\ol{T \xco T'}} $
\xEH
$\xci$
\xEH
$f(X \xcv Y) \xcc f(X) \xcv f(Y)$
\xEH
\xEH
\xEP

$ \xba \xco \xba' \xcn \xbb $
\xEH
\xEH
\xEH
\xEH
\xEH
\xEP

\hline

\xEH
$(PR)$
\xEH
$\xch$
\xEH
$(\xbm PR)$
\xEH
$\xcj$
\xEH
$(eM\xdi)+(I_\xbo)$
\xEP

\cline{3-3}

$ \ol{ \ol{ \xba \xcu \xba' } }$ $ \xcc $ $ \ol{ \ol{ \ol{ \xba } } \xcv
\{ \xba' \}}$
\xEH
$ \ol{ \ol{T \xcv T' } }$ $ \xcc $ $ \ol{ \ol{ \ol{T} } \xcv T' }$
\xEH
$\xci (\xbm dp)+(\xbm\xcc)$
\xEH
$X \xcc Y$ $ \xch $
\xEH
\xEH
\xEP

\cline{3-3}

\xEH
\xEH
$\xcI$ $-(\xbm dp)$
\xEH
$f(Y) \xcs X \xcc f(X)$
\xEH
\xEH
\xEP

\cline{3-3}

\xEH
\xEH
$\xci (\xbm\xcc)$
\xEH
\xEH
\xEH
\xEP

\xEH
\xEH
$T'=\xbf$
\xEH
\xEH
\xEH
\xEP

\cline{3-4}

\xEH
\xEH
$\xci$
\xEH
$(\xbm PR')$
\xEH
\xEH
\xEP

\xEH
\xEH
$T'=\xbf$
\xEH
$f(X) \xcs Y \xcc f(X \xcs Y)$
\xEH
\xEH
\xEP

\hline

$(CUT)$
\xEH
$(CUT)$
\xEH
$\xch$
\xEH
$ (\xbm CUT) $
\xEH
$\xci$
\xEH
$(eM\xdi)+(I_\xbo)$
\xEP

\cline{3-3}
\cline{5-5}

$ T  \xcn \xba; T \xcv \{ \xba\} \xcn \xbb \xch $
\xEH
$T \xcc \ol{T' } \xcc \ol{ \ol{T} }  \xch $
\xEH
$\xci$
\xEH
$f(X) \xcc Y \xcc X  \xch $
\xEH
$\xcH$
\xEH
\xEP

$ T  \xcn \xbb $
\xEH
$ \ol{ \ol{T'} } \xcc \ol{ \ol{T} }$
\xEH
\xEH
$f(X) \xcc f(Y)$
\xEH
\xEH
\xEP

\hline

\end{tabular}
}
\end{turn}
\end{center}
\end{table}

\begin{table}[h]

\index{$(wCM)$}
\index{$(eM\xdf)$}
\index{$(CM_2)$}
\index{$(I_2)$}
\index{$(CM_n)$}
\index{$(I_n)$}
\index{$(CM)$}
\index{Cautious Monotony}
\index{$ (\xbm CM) $}
\index{$(\xdm^+_\xbo)$}
\index{$(ResM)$}
\index{Restricted Monotony}
\index{$(\xbm ResM)$}
\index{$(CUM)$}
\index{Cumulativity}
\index{$(\xbm CUM)$}
\index{$ (\xcc \xcd) $}
\index{$ (\xbm \xcc \xcd) $}
\index{$(eM\xdi)$}
\index{$(I_\xbo)$}
\index{($eM\xdf)$}
\index{$(RatM)$}
\index{Rational Monotony}
\index{$(\xbm RatM)$}
\index{$(\xdm^{++})$}
\index{$(RatM=)$}
\index{$(\xbm =)$}
\index{$(Log=')$}
\index{$(\xbm =')$}
\index{$(DR)$}
\index{$(Log \xFO)$}
\index{$(\xbm \xFO)$}
\index{$(Log \xcv)$}
\index{$(\xbm \xcv)$}
\index{$(Log \xcv')$}
\index{$(\xbm \xcv')$}
\index{$(\xbm \xbe)$}

\tabcolsep=0.5pt
\caption{Logical rules, definitions and connections Part II}

\label{Table Base2-Rules-Def-Conn-2}
\begin{center}
\begin{turn}{90}
{\tiny
\begin{tabular}{|c|c|c|c|c|c|}

\hline

\multicolumn{6}{|c|}{\bf Logical rules, definitions and connections Part II}\\
\hline

\multicolumn{2}{|c|}{Logical rule}
\xEH
Corr.
\xEH
Model set
\xEH
Corr.
\xEH
Size-Rule
\xEP

\hline

\multicolumn{6}{|c|}{Cumulativity}
\xEP

\hline

$(wCM)$
\xEH
\xEH
\xEH
\xEH
trivial
\xEH
$(eM\xdf)$
\xEP

$\xba\xcn\xbb, \xba'\xcl\xba, \xba\xcu\xbb\xcl\xba' \xch $
\xEH
\xEH
\xEH
\xEH
\xEH
\xEP

$\xba'\xcn\xbb$
\xEH
\xEH
\xEH
\xEH
\xEH
\xEP

\hline

$(CM_2)$
\xEH
\xEH
\xEH
\xEH
\xEH
$(I_2)$
\xEP

$\xba\xcn\xbb, \xba\xcn\xbb' \xch \xba\xcu\xbb\xcL\xCN\xbb'$
\xEH
\xEH
\xEH
\xEH
\xEH
\xEP

\hline

$(CM_n)$
\xEH
\xEH
\xEH
\xEH
\xEH
$(I_n)$
\xEP

$\xba\xcn\xbb_1, \ldots, \xba\xcn\xbb_n \xch $
\xEH
\xEH
\xEH
\xEH
\xEH
\xEP

$\xba \xcu \xbb_1 \xcu \ldots \xcu \xbb_{n-1} \xcL\xCN\xbb_n$
\xEH
\xEH
\xEH
\xEH
\xEH
\xEP

\hline

$(CM)$ Cautious Monotony
\xEH
$(CM)$
\xEH
$\xch$
\xEH
$ (\xbm CM) $
\xEH
$\xcj$
\xEH
$(\xdm^+_\xbo)(4)$
\xEP

\cline{3-3}

$ \xba \xcn \xbb,  \xba \xcn \xbb'   \xch $
\xEH
$T \xcc \ol{T' } \xcc \ol{ \ol{T} }  \xch $
\xEH
$\xci$
\xEH
$f(X) \xcc Y \xcc X  \xch $
\xEH
\xEH
\xEP

$ \xba \xcu \xbb \xcn \xbb' $
\xEH
$ \ol{ \ol{T} } \xcc \ol{ \ol{T' } }$
\xEH
\xEH
$f(Y) \xcc f(X)$
\xEH
\xEH
\xEP

\cline{1-1}

\cline{3-4}

or $(ResM)$ Restricted Monotony
\xEH
\xEH
$\xch$
\xEH
$(\xbm ResM)$
\xEH
\xEH
\xEP

\cline{3-3}

$ T  \xcn \xba, \xbb \xch T \xcv \{\xba\} \xcn \xbb $
\xEH
\xEH
$\xci$
\xEH
$ f(X) \xcc A \xcs B \xch $
\xEH
\xEH
\xEP

\xEH
\xEH
\xEH
$f(X \xcs A) \xcc B $
\xEH
\xEH
\xEP

\hline

$(CUM)$ Cumulativity
\xEH
$(CUM)$
\xEH
$\xch$
\xEH
$(\xbm CUM)$
\xEH
$\xci$
\xEH
$(eM\xdi)+(I_\xbo)+(\xdm^{+}_{\xbo})(4)$
\xEP

\cline{3-3}
\cline{5-5}

$ \xba \xcn \xbb   \xch $
\xEH
$T \xcc \ol{T' } \xcc \ol{ \ol{T} }  \xch $
\xEH
$\xci$
\xEH
$f(X) \xcc Y \xcc X  \xch $
\xEH
$\xcH$
\xEH
\xEP

$(\xba \xcn \xbb'   \xcj   \xba \xcu \xbb \xcn \xbb')$
\xEH
$ \ol{ \ol{T} }= \ol{ \ol{T' } }$
\xEH
\xEH
$f(Y)=f(X)$
\xEH
\xEH
\xEP

\hline

\xEH
$ (\xcc \xcd) $
\xEH
$\xch$
\xEH
$ (\xbm \xcc \xcd) $
\xEH
$\xci$
\xEH
$(eM\xdi)+(I_\xbo)+(eM\xdf)$
\xEP

\cline{3-3}
\cline{5-5}

\xEH
$T \xcc \ol{\ol{T'}}, T' \xcc \ol{\ol{T}} \xch $
\xEH
$\xci$
\xEH
$ f(X) \xcc Y, f(Y) \xcc X \xch $
\xEH
$\xcH$
\xEH
\xEP

\xEH
$ \ol{\ol{T'}} = \ol{\ol{T}}$
\xEH
\xEH
$ f(X)=f(Y) $
\xEH
\xEH
\xEP

\hline

\multicolumn{6}{|c|}{Rationality}
\xEP

\hline

$(RatM)$ Rational Monotony
\xEH
$(RatM)$
\xEH
$\xch$
\xEH
$(\xbm RatM)$
\xEH
$\xcj$
\xEH
$(\xdm^{++})$
\xEP

\cline{3-3}

$ \xba \xcn \xbb,  \xba \xcN \xCN \xbb'   \xch $
\xEH
$Con(T \xcv \ol{\ol{T'}})$, $T \xcl T'$ $ \xch $
\xEH
$\xci$ $(\xbm dp)$
\xEH
$X \xcc Y, X \xcs f(Y) \xEd \xCQ   \xch $
\xEH
\xEH
\xEP

\cline{3-3}

$ \xba \xcu \xbb' \xcn \xbb $
\xEH
$ \ol{\ol{T}} \xcd \ol{\ol{\ol{T'}} \xcv T} $
\xEH
$\xcI$ $-(\xbm dp)$
\xEH
$f(X) \xcc f(Y) \xcs X$
\xEH
\xEH
\xEP

\cline{3-3}

\xEH
\xEH
$\xci$ $T=\xbf$
\xEH
\xEH
\xEH
\xEP

\hline

\xEH
$(RatM=)$
\xEH
$\xch$
\xEH
$(\xbm =)$
\xEH
\xEH
\xEP

\cline{3-3}

\xEH
$Con(T \xcv \ol{\ol{T'}})$, $T \xcl T'$ $ \xch $
\xEH
$\xci$ $(\xbm dp)$
\xEH
$X \xcc Y, X \xcs f(Y) \xEd \xCQ   \xch $
\xEH
\xEH
\xEP

\cline{3-3}

\xEH
$ \ol{\ol{T}} = \ol{\ol{\ol{T'}} \xcv T} $
\xEH
$\xcI$ $-(\xbm dp)$
\xEH
$f(X) = f(Y) \xcs X$
\xEH
\xEH
\xEP

\cline{3-3}

\xEH
\xEH
$\xci$ $T=\xbf$
\xEH
\xEH
\xEH
\xEP

\hline

\xEH
$(Log=')$
\xEH
$\xch$
\xEH
$(\xbm =')$
\xEH
\xEH
\xEP

\cline{3-3}

\xEH
$Con(\ol{ \ol{T' } } \xcv T)$ $ \xch $
\xEH
$\xci$ $(\xbm dp)$
\xEH
$f(Y) \xcs X \xEd \xCQ $ $ \xch $
\xEH
\xEH
\xEP

\cline{3-3}

\xEH
$ \ol{ \ol{T \xcv T' } }= \ol{ \ol{ \ol{T' } } \xcv T}$
\xEH
$\xcI$ $-(\xbm dp)$
\xEH
$f(Y \xcs X)=f(Y) \xcs X$
\xEH
\xEH
\xEP

\cline{3-3}

\xEH
\xEH
$\xci$ $T=\xbf$
\xEH
\xEH
\xEH
\xEP

\hline

$(DR)$
\xEH
$(Log \xFO)$
\xEH
$\xch$
\xEH
$(\xbm \xFO)$
\xEH
\xEH
\xEP

\cline{3-3}

$\xba \xco \xbb \xcn \xbg \xch$
\xEH
$ \ol{ \ol{T \xco T' } }$ is one of
\xEH
$\xci$
\xEH
$f(X \xcv Y)$ is one of
\xEH
\xEH
\xEP

$\xba \xcn \xbg$ or $\xbb \xcn \xbg$
\xEH
$\ol{\ol{T}},$ or $\ol{\ol{T'}},$ or $\ol{\ol{T}} \xcs \ol{\ol{T'}}$ (by (CCL))
\xEH
\xEH
$f(X),$ $f(Y)$ or $f(X) \xcv f(Y)$
\xEH
\xEH
\xEP

\hline

\xEH
$(Log \xcv)$
\xEH
$\xch$ $(\xbm\xcc)+(\xbm=)$
\xEH
$(\xbm \xcv)$
\xEH
\xEH
\xEP

\cline{3-3}

\xEH
$Con(\ol{ \ol{T' } } \xcv T),$ $ \xCN Con(\ol{ \ol{T' } }
\xcv \ol{ \ol{T} })$ $ \xch $
\xEH
$\xci$ $(\xbm dp)$
\xEH
$f(Y) \xcs (X-f(X)) \xEd \xCQ $ $ \xch $
\xEH
\xEH
\xEP

\cline{3-3}

\xEH
$ \xCN Con(\ol{ \ol{T \xco T' } } \xcv T')$
\xEH
$\xcI$ $-(\xbm dp)$
\xEH
$f(X \xcv Y) \xcs Y= \xCQ$
\xEH
\xEH
\xEP

\hline

\xEH
$(Log \xcv')$
\xEH
$\xch$ $(\xbm\xcc)+(\xbm=)$
\xEH
$(\xbm \xcv')$
\xEH
\xEH
\xEP

\cline{3-3}

\xEH
$Con(\ol{ \ol{T' } } \xcv T),$ $ \xCN Con(\ol{ \ol{T' }
} \xcv \ol{ \ol{T} })$ $ \xch $
\xEH
$\xci$ $(\xbm dp)$
\xEH
$f(Y) \xcs (X-f(X)) \xEd \xCQ $ $ \xch $
\xEH
\xEH
\xEP

\cline{3-3}

\xEH
$ \ol{ \ol{T \xco T' } }= \ol{ \ol{T} }$
\xEH
$\xcI$ $-(\xbm dp)$
\xEH
$f(X \xcv Y)=f(X)$
\xEH
\xEH
\xEP

\hline

\xEH
\xEH
\xEH
$(\xbm \xbe)$
\xEH
\xEH
\xEP

\xEH
\xEH
\xEH
$a \xbe X-f(X)$ $ \xch $
\xEH
\xEH
\xEP

\xEH
\xEH
\xEH
$ \xcE b \xbe X.a \xce f(\{a,b\})$
\xEH
\xEH
\xEP

\hline

\end{tabular}
}
\end{turn}
\end{center}
\end{table}

$ \xCO $

$ \xCO $
\subsection{
Notation
}

 \xEh

 \xDH

$ \xdp (X)$ is the power set of $X,$ $ \xcc $ is the subset relation, $
\xcb $ the strict part of
$ \xcc,$ i.e. $A \xcb B$ iff $A \xcc B$ and $A \xEd B.$
The operators $ \xcu,$ $ \xCN,$ $ \xco,$ $ \xcp $ and $ \xcl $ have
their usual, classical interpretation.

 \xDH

$ \xdi (X) \xcc \xdp (X)$ and $ \xdf (X) \xcc \xdp (X)$ are dual abstract
notions of size, $ \xdi (X)$ is the
set of ``small'' subsets of $X,$ $ \xdf (X)$ the set of ``big'' subsets of
$X.$ They are
dual in the sense that $A \xbe \xdi (X) \xcj X-A \xbe \xdf (X).$ ``$ \xdi
$'' evokes ``ideal'',
``$ \xdf $'' evokes ``filter'' though the full strength of both is reached
only
in $(< \xbo *s).$ ``s'' evokes ``small'', and ``$(x*s)$'' stands for
``$x$ small sets together are still not everything''.

 \xDH

If $A \xcc X$ is neither in $ \xdi (X),$ nor in $ \xdf (X),$ we say it has
medium size, and
we define $ \xdm (X):= \xdp (X)-(\xdi (X) \xcv \xdf (X)).$ $
\xdm^{+}(X):= \xdp (X)- \xdi (X)$ is the set of subsets
which are not small.

 \xDH

$ \xeA x \xbf $ is a generalized first order quantifier, it is read
``almost all $x$ have property $ \xbf $''. $ \xeA x(\xbf: \xbq)$ is the
relativized version, read:
``almost all $x$ with property $ \xbf $ have also property $ \xbq $''. To
keep the table
``Rules on size''
simple, we write mostly only the non-relativized versions.
Formally, we have $ \xeA x \xbf: \xcj \{x: \xbf (x)\} \xbe \xdf (U)$
where $U$ is the universe, and
$ \xeA x(\xbf: \xbq): \xcj \{x:(\xbf \xcu \xbq)(x)\} \xbe \xdf (\{x:
\xbf (x)\}).$
Soundness and completeness results on $ \xeA $ can be found in
 \cite{Sch95-1}.

 \xDH

Analogously, for propositional logic, we define:

$ \xba \xcn \xbb $ $: \xcj $ $M(\xba \xcu \xbb) \xbe \xdf (M(\xba)),$

where $M(\xbf)$ is the set of models of $ \xbf.$

 \xDH

In preferential structures, $ \xbm (X) \xcc X$ is the set of minimal
elements of $X.$
This generates a principal filter by $ \xdf (X):=\{A \xcc X: \xbm (X) \xcc
A\}.$ Corresponding
properties about $ \xbm $ are not listed systematically.

 \xDH

The usual rules $ \xCf (AND)$ etc. are named here $(AND_{ \xbo }),$ as
they are in a
natural ascending line of similar rules, based on strengthening of the
filter/ideal properties.

 \xDH

For any set of formulas $T,$ and any consequence relation $ \xcn,$ we
will use
$ \ol{T}:=\{ \xbf:T \xcl \xbf \},$ the set of classical consequences of
$T,$ and
$ \ol{ \ol{T} }:=\{ \xbf:T \xcn \xbf \},$ the set of consequences of $T$
under the relation $ \xcn.$

 \xDH

We say that a set $X$ of models is definable by a formula (or a theory)
iff
there is a formula $ \xbf $ (a theory $T)$ such that $X=M(\xbf),$ or
$X=M(T),$ the set of
models of $ \xbf $ or $T,$ respectively.

 \xDH

Most rules are explained in the
table ``Logical rules'',
and ``RW'' stands for Right Weakening.

 \xEj


\subsection{
The groupes of rules
}

The rules concern properties of $ \xdi (X)$ or $ \xdf (X),$ or
dependencies
between such properties for different $X$ and $Y.$ All $X,Y,$ etc. will
be subsets of some universe, say $V.$ Intuitively, $V$ is the set of
all models of some fixed propositional language. It is not
necessary to consider all subsets of $V,$ the intention is to consider
subsets of $V,$ which are definable by a formula or a theory.
So we assume all $X,Y$ etc. taken from some $ \xdy \xcc \xdp (V),$ which
we call the domain.
In the former case, $ \xdy $ is closed under set difference, in
the latter case not necessarily so. (We will mention it when we need some
particular closure property.)

The rules are divided into 5 groups:

 \xEh

 \xDH $ \xCf (Opt),$ which says that ``All'' is optimal - i.e. when there
are no
exceptions, then a soft rule $ \xcn $ holds.

 \xDH 3 monotony rules:

 \xEh
 \xDH $ \xCf (iM)$ is inner monotony, a subset of a small set is small,
 \xDH $(eM \xdi)$ external monotony for ideals: enlarging the base set
keeps small
sets small,
 \xDH $(eM \xdf)$ external monotony for filters: a big subset stays big
when the base
set shrinks.
 \xEj

These three rules are very natural if ``size'' is anything coherent over
change
of base sets. In particular, they can be seen as weakening.

 \xDH $(\xCd)$ keeps proportions, it is here mainly to point the
possibility out.

 \xDH a group of rules $x*s,$ which say how many small sets will not yet
add to
the base set. The notation
``$(< \xbo *s)$'' is an allusion to the full filter property, that
filters are closed under $ \xCf finite$ intersections.
 \xDH Rational monotony, which can best be understood as robustness of $
\xdm^{+},$
see $(\xdm^{++})(3).$

 \xEj

We will assume all base sets to be non-empty in order to avoid pathologies
and
in particular clashes between $ \xCf (Opt)$ and $(1*s).$

Note that the full strength of the usual definitions of
a filter \index{filter}  and an
ideal \index{ideal}  are reached only in line $(< \xbo *s).$
\subsubsection{
Regularities
}

 \xEh

 \xDH

The group of rules $(x*s)$ use ascending strength of $ \xdi / \xdf.$

 \xDH

The column $(\xdm^{+})$ contains interesting algebraic properties. In
particular,
they show a strengthening from $(3*s)$ up to Rationality. They are not
necessarily
equivalent to the corresponding $(I_{x})$ rules, not even in the presence
of the basic rules. The examples show that care has to be taken when
considering
the different variants.

 \xDH

Adding the somewhat superflous $(CM_{2}),$ we have increasing cautious
monotony from $ \xCf (wCM)$ to full $(CM_{ \xbo }).$

 \xDH

We have increasing ``or'' from $ \xCf (wOR)$ to full $(OR_{ \xbo }).$

 \xDH

The line $(2*s)$ is only there because there seems to be no $(
\xdm^{+}_{2}),$ otherwise
we could begin $(n*s)$ at $n=2.$

 \xEj

$ \xCO $
\subsubsection{
Summary
}

We can obtain all rules except $ \xCf (RatM)$ and $(\xCd)$ from $ \xCf
(Opt),$ the monotony
rules - $ \xCf (iM),$ $(eM \xdi),$ $(eM \xdf)$ -, and $(x*s)$ with
increasing $x.$
\subsection{
Table
}

$ \xCO $

The following table is split in two, as it is too big for printing in
one page.

$ \xCO $

\label{Definition Gen-Filter-2-Teile}
\index{Definition Size rules}

(See
Table \ref{Table Base2-Size-Rules-1} (page \pageref{Table Base2-Size-Rules-1}),
``Rules on size - Part $I'' $
and
Table \ref{Table Base2-Size-Rules-2} (page \pageref{Table Base2-Size-Rules-2}),
''Rules on size - Part II".

\begin{table}[h]

\index{Ideal}
\index{Filter}
\index{$ \xdm^+ $}
\index{$ \xeA $}
\index{Optimal proportion}
\index{$(Opt)$}
\index{Monotony}
\index{Improving proportions}
\index{$(iM)$}
\index{internal monotony}
\index{$(eM \xdi)$}
\index{external monotony for ideals}
\index{$(eM \xdf)$}
\index{external monotony for filters}
\index{$(iM)$}
\index{$(eM \xdi)$}
\index{$(eM \xdf)$}
\index{Keeping proportions}
\index{$(\xCd)$}
\index{$(\xdi \xcv disj)$}
\index{$(\xdf \xcv disj)$}
\index{$(\xdm^+ \xcv disj)$}
\index{Robustness of proportions}
\index{$(1*s)$}
\index{$(\xdi_1)$}
\index{$(\xdf_1)$}
\index{$(\xeA_1)$}
\index{$(2*s)$}
\index{$(\xdi_2)$}
\index{$(\xdf_2)$}
\index{$(\xeA_2)$}
\index{$(n*s)$}
\index{$(\xdi_n)$}
\index{$(\xdf_n)$}
\index{$(\xdm^{+}_{n})$}
\index{$(\xeA_n)$}
\index{$(< \xbo*s)$}
\index{$(\xdi_\xbo)$}
\index{$(\xdf_\xbo)$}
\index{$(\xdm^{+}_{ \xbo })$}
\index{$(\xeA_{\xbo})$}
\index{Robustness of $\xdm^+$}
\index{$(\xdm^{++})$}

\caption{Rules on size - Part I}

\label{Table Base2-Size-Rules-1}

\tabcolsep=0.5pt

\begin{center}

{\tiny

\begin{tabular}{|c|c@{.}c|c|c|}

\hline

\multicolumn{5}{|c|}{\bf Rules on size - Part I}\\
\hline

\xEH
``Ideal''
\xEH
``Filter''
\xEH
$ \xdm^+ $
\xEH
$ \xeA $
\xEP

\hline
\hline

\multicolumn{5}{|c|}{Optimal proportion} \xEP

\hline

$(Opt)$
\xEH
$ \xCQ \xbe \xdi (X)$
\xEH
$X \xbe \xdf (X)$
\xEH
\xEH
$ \xcA x \xba \xcp \xeA x \xba$
\xEP

\xEH
\xEH
\xEH
\xEH
\xEP

\hline
\hline

\multicolumn{5}{|c|}
{Monotony (Improving proportions). $(iM)$: internal monotony,}
\xEP
\multicolumn{5}{|c|}
{$(eM \xdi)$: external monotony for ideals,
$(eM \xdf)$: external monotony for filters}
\xEP

\hline

$(iM)$
\xEH
$A \xcc B \xbe \xdi (X)$
\xEH
$A \xbe \xdf (X)$,
\xEH
\xEH
$\xeA x \xba \xcu \xcA x (\xba \xcp \xba')$
\xEP

\xEH
$ \xch $
\xEH
$A \xcc B \xcc X$
\xEH
\xEH
$ \xcp $ $ \xeA x \xba'$
\xEP

\xEH
$A \xbe \xdi (X)$
\xEH
$ \xch $ $B \xbe \xdf (X)$
\xEH
\xEH
\xEP

\hline

$(eM \xdi)$
\xEH
$X \xcc Y \xch$
\xEH
\xEH
\xEH
$\xeA x (\xba: \xbb) \xcu$
\xEP

\xEH
$\xdi (X) \xcc \xdi (Y)$
\xEH
\xEH
\xEH
$\xcA x (\xba' \xcp \xbb) \xcp$
\xEP

\xEH
\xEH
\xEH
\xEH
$\xeA x (\xba \xco \xba': \xbb)$
\xEP

\xEH
\xEH
\xEH
\xEH
\xEP

\xEH
\xEH
\xEH
\xEH
\xEP

\xEH
\xEH
\xEH
\xEH
\xEP

\xEH
\xEH
\xEH
\xEH
\xEP

\hline

$(eM \xdf)$
\xEH
\xEH
$X \xcc Y \xch$
\xEH
\xEH
$\xeA x (\xba: \xbb) \xcu$
\xEP

\xEH
\xEH
$\xdf (Y) \xcs \xdp (X) \xcc $
\xEH
\xEH
$\xcA x (\xbb \xcu \xba \xcp \xba') \xcp$
\xEP

\xEH
\xEH
$ \xdf (X)$
\xEH
\xEH
$\xeA x (\xba \xcu \xba': \xbb)$
\xEP

\xEH
\xEH
\xEH
\xEH
\xEP

\hline
\hline

\multicolumn{5}{|c|}{Keeping proportions} \xEP

\hline

$(\xCd)$
\xEH
$(\xdi \xcv disj)$
\xEH
$(\xdf \xcv disj)$
\xEH
$(\xdm^+ \xcv disj)$
\xEH
$\xeA x(\xba: \xbb) \xcu$
\xEP

\xEH
$A \xbe \xdi (X),$
\xEH
$A \xbe \xdf (X),$
\xEH
$A \xbe \xdm^+ (X),$
\xEH
$\xeA x(\xba': \xbb) \xcu$
\xEP

\xEH
$B \xbe \xdi (Y),$
\xEH
$B \xbe \xdf (Y),$
\xEH
$B \xbe \xdm^+ (Y),$
\xEH
$\xCN \xcE x(\xba \xcu \xba') \xcp$
\xEP

\xEH
$X \xcs Y= \xCQ $ $ \xch $
\xEH
$X \xcs Y= \xCQ $ $ \xch $
\xEH
$X \xcs Y= \xCQ $ $ \xch $
\xEH
$\xeA x(\xba \xco \xba': \xbb)$
\xEP

\xEH
$A \xcv B \xbe \xdi (X \xcv Y)$
\xEH
$A \xcv B \xbe \xdf (X \xcv Y)$
\xEH
$A \xcv B \xbe \xdm^+ (X \xcv Y)$
\xEH
\xEP

\xEH
\xEH
\xEH
\xEH
\xEP

\xEH
\xEH
\xEH
\xEH
\xEP

\hline
\hline

\multicolumn{5}{|c|}{Robustness of proportions: $n*small \xEd All$} \xEP

\hline

$(1*s)$
\xEH
$(\xdi_1)$
\xEH
$(\xdf_1)$
\xEH
\xEH
$(\xeA_1)$
\xEP

\xEH
$X \xce \xdi (X)$
\xEH
$ \xCQ \xce \xdf (X)$
\xEH
\xEH
$ \xeA x \xba \xcp \xcE x \xba $
\xEP

\hline

$(2*s)$
\xEH
$(\xdi_2)$
\xEH
$(\xdf_2)$
\xEH
\xEH
$(\xeA_2)$
\xEP

\xEH
$A,B \xbe \xdi (X) \xch $
\xEH
$A,B \xbe \xdf (X) \xch $
\xEH
\xEH
$ \xeA x \xba \xcu \xeA x \xbb $
\xEP

\xEH
$ A \xcv B \xEd X$
\xEH
$A \xcs B \xEd \xCQ $
\xEH
\xEH
$ \xcp $ $ \xcE x(\xba \xcu \xbb)$
\xEP

\hline

$(n*s)$
\xEH
$(\xdi_n)$
\xEH
$(\xdf_n)$
\xEH
$(\xdm^{+}_{n})$
\xEH
$(\xeA_n)$
\xEP

$(n \xcg 3)$
\xEH
$A_{1},.,A_{n} \xbe \xdi (X) $
\xEH
$A_{1},.,A_{n} \xbe \xdi (X) $
\xEH
$X_{1} \xbe \xdf (X_{2}),., $
\xEH
$ \xeA x \xba_{1} \xcu.  \xcu \xeA x \xba_{n} $
\xEP

\xEH
$ \xch $
\xEH
$ \xch $
\xEH
$ X_{n-1} \xbe \xdf (X_{n})$ $ \xch $
\xEH
$ \xcp $
\xEP

\xEH
$ A_{1} \xcv.  \xcv A_{n} \xEd X $
\xEH
$A_{1} \xcs.  \xcs A_{n} \xEd \xCQ$
\xEH
$X_{1} \xbe \xdm^{+}(X_{n})$
\xEH
$ \xcE x (\xba_{1} \xcu.  \xcu \xba_{n}) $
\xEP

\xEH
\xEH
\xEH
\xEH
\xEP

\hline

$(< \xbo*s)$
\xEH
$(\xdi_\xbo)$
\xEH
$(\xdf_\xbo)$
\xEH
$(\xdm^{+}_{ \xbo })$
\xEH
$(\xeA_{\xbo})$
\xEP

\xEH
$A,B \xbe \xdi (X) \xch $
\xEH
$A,B \xbe \xdf (X) \xch $
\xEH
(1)
\xEH
$ \xeA x \xba \xcu \xeA x \xbb \xcp $
\xEP

\xEH
$ A \xcv B \xbe \xdi (X)$
\xEH
$ A \xcs B \xbe \xdf (X)$
\xEH
$A \xbe \xdf (X),$ $X \xbe \xdm^{+}(Y)$
\xEH
$ \xeA x(\xba \xcu \xbb)$
\xEP

\xEH
\xEH
\xEH
$ \xch $ $A \xbe \xdm^{+}(Y)$
\xEH
\xEP

\xEH
\xEH
\xEH
(2)
\xEH
\xEP

\xEH
\xEH
\xEH
$A \xbe \xdm^{+}(X),$ $X \xbe \xdf (Y)$
\xEH
\xEP

\xEH
\xEH
\xEH
$ \xch $ $A \xbe \xdm^{+}(Y)$
\xEH
\xEP

\xEH
\xEH
\xEH
(3)
\xEH
\xEP

\xEH
\xEH
\xEH
$A \xbe \xdf (X),$ $X \xbe \xdf (Y)$
\xEH
\xEP

\xEH
\xEH
\xEH
$ \xch $ $A \xbe \xdf (Y)$
\xEH
\xEP

\xEH
\xEH
\xEH
(4)
\xEH
\xEP

\xEH
\xEH
\xEH
$A,B \xbe \xdi (X)$ $ \xch $
\xEH
\xEP

\xEH
\xEH
\xEH
$A-B \xbe \xdi (X-$B)
\xEH
\xEP

\hline
\hline

\multicolumn{5}{|c|}{Robustness of $\xdm^+$} \xEP

\hline

$(\xdm^{++})$
\xEH
\xEH
\xEH
$(\xdm^{++})$
\xEH
\xEP

\xEH
\xEH
\xEH
(1)
\xEH
\xEP

\xEH
\xEH
\xEH
$A \xbe \xdi (X),$ $B \xce \xdf (X)$
\xEH
\xEP

\xEH
\xEH
\xEH
$ \xch $ $A-B \xbe \xdi (X-B)$
\xEH
\xEP

\xEH
\xEH
\xEH
(2)
\xEH
\xEP

\xEH
\xEH
\xEH
$A \xbe \xdf (X), B \xce \xdf (X)$
\xEH
\xEP

\xEH
\xEH
\xEH
$ \xch $ $A-B \xbe \xdf (X-B)$
\xEH
\xEP

\xEH
\xEH
\xEH
(3)
\xEH
\xEP

\xEH
\xEH
\xEH
$A \xbe \xdm^+ (X),$
\xEH
\xEP

\xEH
\xEH
\xEH
$X \xbe \xdm^+ (Y)$
\xEH
\xEP

\xEH
\xEH
\xEH
$ \xch $ $A \xbe \xdm^+ (Y)$
\xEH
\xEP

\hline

\end{tabular}
}
\end{center}
\end{table}

\newpage

\begin{table}[h]

\index{AND}
\index{OR}
\index{Cautious Monotony}
\index{Rational Monotony}
\index{Optimal proportion}
\index{$(Opt)$}
\index{$(SC)$}
\index{Monotony}
\index{Improving proportions}
\index{$(iM)$}
\index{$(RW)$}
\index{$(eM \xdi)$}
\index{$(PR')$}
\index{$(wOR)$}
\index{$(\xbm wOR)$}
\index{$(\xbm PR)$}
\index{$(eM \xdf)$}
\index{$(wCM)$}
\index{Keeping proportions}
\index{$(\xCd)$}
\index{$(NR)$}
\index{$(disjOR)$}
\index{$(\xbm disjOR)$}
\index{Robustness of proportions}
\index{$(1*s)$}
\index{$(CP)$}
\index{$(AND_{1})$}
\index{$(2*s)$}
\index{$(AND_{2})$}
\index{$(OR_{2})$}
\index{$(CM_{2})$}
\index{$(n*s)$}
\index{$(AND_{n})$}
\index{$(OR_{n})$}
\index{$(CM_{n})$}
\index{$(< \xbo*s)$}
\index{$(AND_{ \xbo })$}
\index{$(OR_{ \xbo })$}
\index{$(CM_{ \xbo })$}
\index{$(\xbm OR)$}
\index{$(\xbm CM)$}
\index{Robustness of $\xdm^+$}
\index{$(\xdm^{++})$}
\index{$(RatM)$}
\index{$(\xbm RatM)$}

\caption{Rules on size - Part II}

\label{Table Base2-Size-Rules-2}

\begin{center}

{\tiny

\begin{tabular}{|c|c|c|c|c|}

\hline

\multicolumn{5}{|c|}{\bf Rules on size - Part II}\\
\hline

\xEH
various rules
\xEH
AND
\xEH
OR
\xEH
Caut./Rat.Mon.
\xEP

\hline
\hline

\multicolumn{5}{|c|}{Optimal proportion} \xEP

\hline

$(Opt)$
\xEH
$(SC)$
\xEH
\xEH
\xEH
\xEP

\xEH
$ \xba \xcl \xbb \xch \xba \xcn \xbb $
\xEH
\xEH
\xEH
\xEP

\hline
\hline

\multicolumn{5}{|c|}
{Monotony (Improving proportions)}
\xEP
\multicolumn{5}{|c|}
{}
\xEP

\hline

$(iM)$
\xEH
$(RW)$
\xEH
\xEH
\xEH
\xEP

\xEH
$ \xba \xcn \xbb, \xbb \xcl \xbb' \xch $
\xEH
\xEH
\xEH
\xEP

\xEH
$ \xba \xcn \xbb' $
\xEH
\xEH
\xEH
\xEP

\hline

$(eM \xdi)$
\xEH
$(PR')$
\xEH
\xEH
$(wOR)$
\xEH
\xEP

\xEH
$\xba \xcn \xbb, \xba \xcl \xba',$
\xEH
\xEH
$ \xba \xcn \xbb,$ $ \xba' \xcl \xbb $ $ \xch $
\xEH
\xEP

\xEH
$\xba' \xcu \xCN \xba \xcl \xbb \xch$
\xEH
\xEH
$ \xba \xco \xba' \xcn \xbb $
\xEH
\xEP

\xEH
$\xba' \xcn \xbb$
\xEH
\xEH
$(\xbm wOR)$
\xEH
\xEP

\xEH
$(\xbm PR)$
\xEH
\xEH
$\xbm(X \xcv Y) \xcc \xbm(X) \xcv Y$
\xEH
\xEP

\xEH
$X \xcc Y \xch$
\xEH
\xEH
\xEH
\xEP

\xEH
$\xbm(Y) \xcs X \xcc \xbm(X)$
\xEH
\xEH
\xEH
\xEP

\hline

$(eM \xdf)$
\xEH
\xEH
\xEH
\xEH
$(wCM)$
\xEP

\xEH
\xEH
\xEH
\xEH
$\xba \xcn \xbb, \xba' \xcl \xba,$
\xEP

\xEH
\xEH
\xEH
\xEH
$\xba \xcu \xbb \xcl \xba' \xch$
\xEP

\xEH
\xEH
\xEH
\xEH
$\xba' \xcn \xbb$
\xEP

\hline
\hline

\multicolumn{5}{|c|}{Keeping proportions} \xEP

\hline

$(\xCd)$
\xEH
$(NR)$
\xEH
\xEH
$(disjOR)$
\xEH
\xEP

\xEH
$\xba \xcn \xbb \xch$
\xEH
\xEH
$ \xba \xcn \xbb,$ $ \xba' \xcn \xbb' $
\xEH
\xEP

\xEH
$\xba \xcu \xbg \xcn \xbb$
\xEH
\xEH
$ \xba \xcl \xCN \xba',$ $ \xch $
\xEH
\xEP

\xEH
or
\xEH
\xEH
$ \xba \xco \xba' \xcn \xbb \xco \xbb' $
\xEH
\xEP

\xEH
$\xba \xcu \xCN \xbg \xcn \xbb$
\xEH
\xEH
$(\xbm disjOR)$
\xEH
\xEP

\xEH
\xEH
\xEH
$X \xcs Y = \xCQ \xch$
\xEH
\xEP

\xEH
\xEH
\xEH
$\xbm(X \xcv Y) \xcc \xbm(X) \xcv \xbm(Y)$
\xEH
\xEP

\hline
\hline

\multicolumn{5}{|c|}{Robustness of proportions: $n*small \xEd All$} \xEP

\hline

$(1*s)$
\xEH
$(CP)$
\xEH
$(AND_{1})$
\xEH
\xEH
\xEP

\xEH
$\xba\xcn\xcT \xch \xba\xcl\xcT$
\xEH
$ \xba \xcn \xbb $ $ \xch $ $ \xba \xcL \xCN \xbb $
\xEH
\xEH
\xEP

\hline

$(2*s)$
\xEH
\xEH
$(AND_{2})$
\xEH
$(OR_{2})$
\xEH
$(CM_{2})$
\xEP

\xEH
\xEH
$ \xba \xcn \xbb,$ $ \xba \xcn \xbb' $ $ \xch $
\xEH
$ \xba \xcn \xbb \xch \xba \xcN \xCN \xbb $
\xEH
$ \xba \xcn \xbb \xch \xba \xcN \xCN \xbb $
\xEP

\xEH
\xEH
$ \xba \xcL \xCN \xbb \xco \xCN \xbb' $
\xEH
\xEH
\xEP

\hline

$(n*s)$
\xEH
\xEH
$(AND_{n})$
\xEH
$(OR_{n})$
\xEH
$(CM_{n})$
\xEP

$(n \xcg 3)$
\xEH
\xEH
$ \xba \xcn \xbb_{1},., \xba \xcn \xbb_{n}$
\xEH
$ \xba_{1} \xcn \xbb,., \xba_{n-1} \xcn \xbb $
\xEH
$ \xba \xcn \xbb_{1},., \xba \xcn \xbb_{n-1}$
\xEP

\xEH
\xEH
$ \xch $
\xEH
$ \xch $
\xEH
$ \xch $
\xEP

\xEH
\xEH
$ \xba \xcL \xCN \xbb_{1} \xco.  \xco \xCN \xbb_{n}$
\xEH
$ \xba_{1} \xco.  \xco \xba_{n-1} \xcN \xCN \xbb $
\xEH
$ \xba \xcu \xbb_1 \xcu.  \xcu \xbb_{n-2} \xcN $
\xEP

\xEH
\xEH
\xEH
\xEH
$ \xCN \xbb_{n-1}$
\xEP

\hline

$(< \xbo*s)$
\xEH
\xEH
$(AND_{ \xbo })$
\xEH
$(OR_{ \xbo })$
\xEH
$(CM_{ \xbo })$
\xEP

\xEH
\xEH
$ \xba \xcn \xbb,$ $ \xba \xcn \xbb' $ $ \xch $
\xEH
$ \xba \xcn \xbb,$ $ \xba' \xcn \xbb $ $ \xch $
\xEH
$ \xba \xcn \xbb,$ $ \xba \xcn \xbb' $ $ \xch $
\xEP

\xEH
\xEH
$ \xba \xcn \xbb \xcu \xbb' $
\xEH
$ \xba \xco \xba' \xcn \xbb $
\xEH
$ \xba \xcu \xbb \xcn \xbb' $
\xEP

\xEH
\xEH
\xEH
$(\xbm OR)$
\xEH
$(\xbm CM)$
\xEP

\xEH
\xEH
\xEH
$\xbm(X \xcv Y) \xcc \xbm(X) \xcv \xbm(Y)$
\xEH
$\xbm(X) \xcc Y \xcc X \xch$
\xEP

\xEH
\xEH
\xEH
\xEH
$\xbm(Y) \xcc \xbm(X)$
\xEP

\xEH
\xEH
\xEH
\xEH
\xEP

\xEH
\xEH
\xEH
\xEH
\xEP

\xEH
\xEH
\xEH
\xEH
\xEP

\xEH
\xEH
\xEH
\xEH
\xEP

\xEH
\xEH
\xEH
\xEH
\xEP

\xEH
\xEH
\xEH
\xEH
\xEP

\xEH
\xEH
\xEH
\xEH
\xEP

\hline
\hline

\multicolumn{5}{|c|}{Robustness of $\xdm^+$} \xEP

\hline

$(\xdm^{++})$
\xEH
\xEH
\xEH
\xEH
$(RatM)$
\xEP

\xEH
\xEH
\xEH
\xEH
$ \xba \xcn \xbb,  \xba \xcN \xCN \xbb'   \xch $
\xEP

\xEH
\xEH
\xEH
\xEH
$ \xba \xcu \xbb' \xcn \xbb $
\xEP

\xEH
\xEH
\xEH
\xEH
$(\xbm RatM)$
\xEP

\xEH
\xEH
\xEH
\xEH
$X \xcc Y,$
\xEP

\xEH
\xEH
\xEH
\xEH
$X \xcs \xbm(Y) \xEd \xCQ \xch$
\xEP

\xEH
\xEH
\xEH
\xEH
$\xbm(X) \xcc \xbm(Y) \xcs X$
\xEP

\xEH
\xEH
\xEH
\xEH
\xEP

\xEH
\xEH
\xEH
\xEH
\xEP

\xEH
\xEH
\xEH
\xEH
\xEP

\xEH
\xEH
\xEH
\xEH
\xEP

\hline

\end{tabular}
}
\end{center}
\end{table}

$ \xCO $

$ \xCO $
\clearpage
\section{
Multiplication of size for subsets
}

Here we have nested sets, $A \xcc X \xcc Y,$ $ \xCf A$ is a certain
proportion of $X,$ and $X$ of $Y,$
resulting in a multiplication of relative size or proportions. This is a
classical subject of nonmonotonic logic,
see the last section, taken from  \cite{GS09a}, it is partly
repeated here
to stress the common points with the other scenario.

$ \xCO $

\vspace{10mm}

\begin{diagram}

\label{Diagram Mul-Add}
\index{Diagram Mul-Add}

\centering
\setlength{\unitlength}{1mm}
{\renewcommand{\dashlinestretch}{30}
\begin{picture}(150,100)(0,0)

\path(20,90)(90,90)(90,30)(20,30)(20,90)

\path(55,90)(55,30)
\path(55,60)(90,60)

\path(55,20)(90,20)
\path(55,22)(55,18)
\path(90,22)(90,18)

\path(20,10)(90,10)
\path(20,12)(20,8)
\path(90,12)(90,8)

\path(100,60)(100,30)
\path(98,60)(102,60)
\path(98,30)(102,30)

\put(54,6){{\xssc $Y$}}
\put(72,16){{\xssc $X$}}
\put(102,44){{\xssc $A$}}

\put(50,1){{\xssc Scenario 1}}

\end{picture}
}

\end{diagram}

\vspace{4mm}

$ \xCO $
\subsection{
Properties
}

Diagram \ref{Diagram Mul-Add} (page \pageref{Diagram Mul-Add})  is to be read as
follows:
The whole set $Y$ is split in $X$ and $Y- \xCf X,$ $X$ is split in $ \xCf
A$ and $X- \xCf A.$
$X$ is a small/medium/big part of $Y,$
$ \xCf A$ is a small/medium/big part of $X.$
The question is: is $ \xCf A$ a small/medium/big part of $Y?$

Note that the relation of $ \xCf A$ to $X$ is conceptually different from
that of
$X$ to $ \xCf Y,$ as we change the base set by going from $X$ to $ \xCf
Y,$ but not when going
from $ \xCf A$ to $X.$ Thus, in particular, when we read the diagram as
expressing
multiplication, commutativity is not necessarily true.

We looked at this scenario already in  \cite{GS09a}, but there from
an
additive point of view, using various basic properties like (iM), $(eM
\xdi),$
$(eM \xdf).$ Here, we use just multiplication - except sometimes for
motivation.

We examine different rules:

If $Y=X$ or $X=A,$ there is nothing to show, so 1 is the neutral element
of
multiplication.

If $X \xbe \xdi (Y)$ or $A \xbe \xdi (X),$ then we should have $A \xbe
\xdi (Y).$ (Use for motivation
(iM) or $(eM \xdi)$ respectively.)

So it remains to look at the following cases, with the ``natural'' answers
given already:

(1) $X \xbe \xdf (Y),$ $A \xbe \xdf (X)$ $ \xch $ $A \xbe \xdf (Y),$

(2) $X \xbe \xdm^{+}(Y),$ $A \xbe \xdf (X)$ $ \xch $ $A \xbe \xdm^{+}(Y),$

(3) $X \xbe \xdf (Y),$ $A \xbe \xdm^{+}(X)$ $ \xch $ $A \xbe \xdm^{+}(Y),$

(4) $X \xbe \xdm^{+}(Y),$ $A \xbe \xdm^{+}(X)$ $ \xch $ $A \xbe
\xdm^{+}(Y).$

But (1) is case (3) of $(\xdm^{+}_{ \xbo })$ in  \cite{GS09a}, see
Table ``Rules on size'' in
Section \ref{Section Mul-Intro} (page \pageref{Section Mul-Intro}).

(2) is case (2) of $(\xdm^{+}_{ \xbo })$ there,

(3) is case (1) of $(\xdm^{+}_{ \xbo })$ there, finally,

(4) is $(\xdm^{++})$ there.

So the first three correspond to various expressions of $(AND_{ \xbo }),$
$(OR_{ \xbo }),$
$(CM_{ \xbo }),$ the last one to $ \xCf (RatM).$

But we can read them also the other way round, e.g.:

(1) corresponds to: $ \xba \xcn \xbb,$ $ \xba \xcu \xbb \xcn \xbg $ $
\xch $ $ \xba \xcn \xbg,$

(2) corresponds to: $ \xba \xcN \xCN \xbb,$ $ \xba \xcu \xbb \xcn \xbg $
$ \xch $ $ \xba \xcN \xCN (\xbb \xcu \xbg),$

(3) corresponds to: $ \xba \xcn \xbb,$ $ \xba \xcu \xbb \xcN \xCN \xbg $
$ \xch $ $ \xba \xcN \xCN (\xbb \xcu \xbg).$

All these rules might be seen as too idealistic, so just as we did in
 \cite{GS09a}, we can consider milder versions:
We might for instance consider a rule which says that $big* \Xl *big,$ $n$
times,
is not small. Consider for instance the case $n=2.$
So we would conclude that $ \xCf A$ is not small in $Y.$ In terms of
logic, we
then have: $ \xba \xcn \xbb,$ $ \xba \xcu \xbb \xcn \xbg $ $ \xch $ $
\xba \xcN (\xCN \xbb \xco \xCN \xbg).$ We can obtain the same
logical property from $3*small \xEd all.$
\section{
Multiplication of size for subspaces
}
\subsection{
Properties
}

\label{Section Mul-Mul}

$ \xCO $

\vspace{10mm}

\begin{diagram}

\label{Diagram Mul-Prod}
\index{Diagram Mul-Prod}

\centering
\setlength{\unitlength}{1mm}
{\renewcommand{\dashlinestretch}{30}
\begin{picture}(150,100)(0,0)

\path(20,90)(90,90)(90,30)(20,30)(20,90)

\path(55,90)(55,30)

\path(55,70)(90,70)
\path(20,75)(55,75)

\path(20,20)(90,20)

\path(20,22)(20,18)
\path(55,22)(55,18)
\path(90,22)(90,18)

\path(10,75)(10,30)
\path(8,75)(12,75)
\path(8,30)(12,30)

\path(100,70)(100,30)
\path(98,70)(102,70)
\path(98,30)(102,30)

\put(37,16){{\xssc $\xbS_1$}}
\put(72,16){{\xssc $\xbS_2$}}

\put(5,50){{\xssc $\xbG_1$}}
\put(101,50){{\xssc $\xbG_2$}}

\put(50,1){{\xssc Scenario 2}}

\end{picture}
}

\end{diagram}

\vspace{4mm}

$ \xCO $

In this scenario, $ \xbS_{i}$ are sets of sequences,
see Diagram \ref{Diagram Mul-Prod} (page \pageref{Diagram Mul-Prod}).
(Correponding, intuitively, to a set of models in language $ \xdl_{i},$
$ \xbS_{i}$ will be the set of $ \xba_{i}-$models, and the subsets $
\xbG_{i}$ are to be seen as
the ``best'' models, where $ \xbb_{i}$ will hold. The languages are supposed
to
be disjoint sublanguages of a common language $ \xdl.)$

In this scenario, the $ \xbS_{i}$ have symmetrical roles,
so there is no intuitive reason for multiplication not to be commutative.

We can interpret the situation twofold:

First, we work separately in sublanguage $ \xdl_{1}$ and $ \xdl_{2},$ and,
say, $ \xba_{i}$ and $ \xbb_{i}$
are both defined in $ \xdl_{i},$ and we look at $ \xba_{i} \xcn \xbb_{i}$
in the sublanguage $ \xdl_{i},$
or, we consider both $ \xba_{i}$ and $ \xbb_{i}$ in the big language $
\xdl,$ and look at
$ \xba_{i} \xcn \xbb_{i}$ in $ \xdl.$ These two ways are a priori
completely different.
Speaking in preferential terms, it is not at all clear why the orderings
on the submodels should have anything to do with the orderings on the
whole models. It seems a very desirable property, but we have to postulate
it, which we do now (an overview is given
in Table \ref{Table Mul-Laws} (page \pageref{Table Mul-Laws})). We give now
informally a list of such
rules,
mainly to show the connection with the first scenario. Later,
see Definition \ref{Definition Mul-Size-Rules} (page \pageref{Definition
Mul-Size-Rules}), we will
introduce formally some rules for which we show a connection with
interpolation.
Here, e.g.,
``$(big*big \xch big)$'' stands for
``if both factors are big, so will be the product'',
this will be abbreviated by
``$b*b \xch b$'' in Table \ref{Table Mul-Laws} (page \pageref{Table Mul-Laws})
.

$(big*1 \xch big)$ Let $ \xbG_{1} \xcc \xbS_{1},$ if $ \xbG_{1} \xbe \xdf
(\xbS_{1}),$ then $ \xbG_{1} \xCK \xbS_{2} \xbe \xdf (\xbS_{1} \xCK
\xbS_{2}),$
(and the dual rule for $ \xbS_{2}$ and $ \xbG_{2}).$

This property preserves proportions, so it seems intuitively quite
uncontested, whenever we admit coherence over products. (Recall that there
was nothing to show in the first scenario.)

When we re-consider above case: suppose $ \xba \xcn \xbb $ in the
sublanguage, so
$M(\xbb) \xbe \xdf (M(\xba))$ in the sublanguage, so by $(big*1 \xch
big),$ $M(\xbb) \xbe \xdf (M(\xba))$
in the big language $ \xdl.$

We obtain the dual rule for small (and likewise, medium size) sets:

$(small*1 \xch small)$ Let $ \xbG_{1} \xcc \xbS_{1},$ if $ \xbG_{1} \xbe
\xdi (\xbS_{1}),$ then $ \xbG_{1} \xCK \xbS_{2} \xbe \xdi (\xbS_{1} \xCK
\xbS_{2}),$
(and the dual rule for $ \xbS_{2}$ and $ \xbG_{2}),$

establishing $ \xCf All=1$ as the neutral element for multiplication.

We look now at other, plausible rules:

$(small*x \xch small)$ $ \xbG_{1} \xbe \xdi (\xbS_{1}),$ $ \xbG_{2} \xcc
\xbS_{2}$ $ \xch $ $ \xbG_{1} \xCK \xbG_{2} \xbe \xdi (\xbS_{1} \xCK
\xbS_{2})$

$(big*big \xch big)$ $ \xbG_{1} \xbe \xdf (\xbS_{1}),$ $ \xbG_{2} \xbe
\xdf (\xbS_{2})$ $ \xch $ $ \xbG_{1} \xCK \xbG_{2} \xbe \xdf (\xbS_{1}
\xCK \xbS_{2})$

$(big*medium \xch medium)$ $ \xbG_{1} \xbe \xdf (\xbS_{1}),$ $ \xbG_{2}
\xbe \xdm^{+}(\xbS_{2})$ $ \xch $ $ \xbG_{1} \xCK \xbG_{2} \xbe \xdm^{+}(
\xbS_{1} \xCK \xbS_{2})$

$(medium*medium \xch medium)$ $ \xbG_{1} \xbe \xdm^{+}(\xbS_{1}),$ $
\xbG_{2} \xbe \xdm^{+}(\xbS_{2})$ $ \xch $ $ \xbG_{1} \xCK \xbG_{2} \xbe
\xdm^{+}(\xbS_{1} \xCK \xbS_{2})$

When we accept all above rules, we can invert $(big*big \xch big),$ as a
big product
must be composed of big components. Likewise, at least one component of a
small
product has to be small - see
Fact \ref{Fact Mul-Big-Small} (page \pageref{Fact Mul-Big-Small}).

We see that these properties give a lot of modularity. We can calculate
the consequences of $ \xba $ and $ \xba' $ separately - provided $ \xba
,$ $ \xba' $ use disjoint
alphabets - and put the results together afterwards. Such properties are
particularly interesting for classification purposes, where subclasses
are defined with disjoint alphabets.
\subsection{
Size multiplication and corresponding preferential relations
}

We turn to those conditions which provide
the key to non-monotonic interpolation theorems - see
Section \ref{Section Mul-Int} (page \pageref{Section Mul-Int}).
We quote from  \cite{GS09c} the following pairwise equivalent
conditions
$(S*1),$ $(\xbm *1),$ $(S*2),$ $(\xbm *2),$ and add a new condition,
$(s*s),$
for a principal filter generated by $ \xbm:$

\bd

$\hspace{0.01em}$


\label{Definition Mul-Size-Rules}

$(S*1)$ $ \xbD \xcc \xbS' \xCK \xbS'' $ is big iff there is $ \xbG =
\xbG' \xCK \xbG'' \xcc \xbD $ s.t. $ \xbG' \xcc \xbS' $ and
$ \xbG'' \xcc \xbS'' $ are big

$(\xbm *1)$ $ \xbm (\xbS' \xCK \xbS'')= \xbm (\xbS') \xCK \xbm (
\xbS'')$

$(S*2)$ $ \xbG \xcc \xbS $ is big $ \xch $ $ \xbG \xex X' \xcc \xbS \xex
X' $ is big - where
$ \xbS $ is not necessarily a product.

$(\xbm *2)$ $ \xbm (\xbS) \xcc \xbG $ $ \xch $ $ \xbm (\xbS \xex X')
\xcc \xbG \xex X' $

$(s*s)$ Let $ \xbG_{i} \xcc \xbS_{i},$ then $ \xbG_{1} \xCK \xbG_{2} \xcc
\xbS_{1} \xCK \xbS_{2}$ is small iff $ \xbG_{1} \xcc \xbS_{1}$ is small or
$ \xbG_{1} \xcc \xbS_{1}$ is small.

\ed

$(\xbm *1)$ and $(s*s)$ are equivalent in the following sense:

\bfa

$\hspace{0.01em}$


\label{Fact Mul-Big-Small}

Let the notion of size satisfy $ \xCf (Opt),$ $ \xCf (iM),$ and $(< \xbo
*s),$ see
the tables ``Rules on size'' in
Section \ref{Section Mul-Intro} (page \pageref{Section Mul-Intro}).
Then $(\xbm *1)$ and $(s*s)$ are equivalent.

\efa

\subparagraph{
Proof
}

$\hspace{0.01em}$


``$ \xch $'':

(1) Let $ \xbG' \xcc \xbS' $ be small, we show that $ \xbG' \xCK \xbG
'' \xcc \xbS' \xCK \xbS'' $ is small.
So $ \xbS' - \xbG' \xcc \xbS' $ is big, so by $ \xCf (Opt)$ and $(\xbm
*1)$
$(\xbS' - \xbG') \xCK \xbS'' \xcc \xbS' \xCK \xbS'' $ is big, so
$ \xbG' \xCK \xbS'' $ $=$ $(\xbS' \xCK \xbS'')-((\xbS' - \xbG')
\xCK \xbS'')$ $ \xcc $ $ \xbS' \xCK \xbS'' $ is small,
so by $ \xCf (iM)$ $ \xbG' \xCK \xbG'' \xcc \xbS' \xCK \xbS'' $ is
small.

(2) Suppose $ \xbG' \xcc \xbS' $ and $ \xbG'' \xcc \xbS'' $ are not
small, we show that $ \xbG' \xCK \xbG'' \xcc \xbS' \xCK \xbS'' $
is not small. So $ \xbS' - \xbG' \xcc \xbS' $ and $ \xbS'' - \xbG''
\xcc \xbS'' $ are not big.
We show that $Z$ $:=$ $((\xbS' \xCK \xbS'')-(\xbG' \xCK \xbG''))$
$ \xcc $ $ \xbS' \xCK \xbS'' $ is not big.
$Z$ $=$ $(\xbS' \xCK (\xbS'' - \xbG'')) \xcv ((\xbS' - \xbG')
\xCK \xbS'').$

Suppose $X' \xCK X'' $ $ \xcc $ $Z,$ then $X' \xcc \xbS' - \xbG' $ or
$X'' \xcc \xbS'' - \xbG''.$ Proof:
Let $X' \xcC \xbS' - \xbG' $ and $X'' \xcC \xbS'' - \xbG'',$ but $X'
\xCK X'' \xcc Z.$ Let $ \xbs' \xbe X' -(\xbS' - \xbG'),$
$ \xbs'' \xbe X'' -(\xbS'' - \xbG''),$ consider $ \xbs' \xbs''.$
$ \xbs' \xbs'' \xce (\xbS' - \xbG') \xCK \xbS'',$ as $ \xbs' \xce
\xbS' - \xbG',$
$ \xbs' \xbs'' \xce \xbS' \xCK (\xbS'' \xCK \xbG''),$ as $ \xbs''
\xce \xbS'' - \xbG'',$ so $ \xbs' \xbs'' \xce Z.$

By prerequisite, $ \xbS' - \xbG' \xcc \xbS' $ is not big, $ \xbS'' -
\xbG'' \xcc \xbS'' $ is not big, so by
$ \xCf (iM)$ no $X' $ with $X' \xcc \xbS' - \xbG' $ is big, no $X'' $
with $X'' \xcc \xbS'' - \xbG'' $ is big, so
by $(\xbm *1)$ or $(S*1)$ $Z \xcc \xbS' \xCK \xbS'' $ is not big, so $
\xbG' \xCK \xbG'' \xcc \xbS' \xCK \xbS'' $
is not small.

``$ \xci $'':

(1) Suppose $ \xbG' \xcc \xbS' $ is big, $ \xbG'' \xcc \xbS'' $ is
big, we have to show
$ \xbG' \xCK \xbG'' \xcc \xbS' \xCK \xbS'' $ is big.
$ \xbS' - \xbG' \xcc \xbS' $ is small, $ \xbS'' - \xbG'' \xcc \xbS''
$ is small, so by $(s*s)$
$(\xbS' - \xbG') \xCK \xbS'' \xcc \xbS' \xCK \xbS'' $ is small and
$ \xbS' \xCK (\xbS'' - \xbG'') \xcc \xbS' \xCK \xbS'' $ is small,
so by $(< \xbo *s)$
$(\xbS' \xCK \xbS'')-(\xbG' \xCK \xbG'')$ $=$
$((\xbS' - \xbG') \xCK \xbS'') \xcv (\xbS' \xCK (\xbS'' - \xbG
''))$ $ \xcc $ $ \xbS' \xCK \xbS'' $ is small,
so $ \xbG' \xCK \xbG'' \xcc \xbS' \xCK \xbS'' $ is big.

(2) Suppose $ \xbG' \xCK \xbG'' \xcc \xbS' \xCK \xbS'' $ is big, we
have to show $ \xbG' \xcc \xbS' $ is big, and
$ \xbG'' \xcc \xbS'' $ is big. By prerequisite,
$(\xbS' \xCK \xbS'')-(\xbG' \xCK \xbG'')$ $=$
$((\xbS' - \xbG') \xCK \xbS'') \xcv (\xbS' \xCK (\xbS'' - \xbG
''))$ $ \xcc $ $ \xbS' \xCK \xbS'' $ is small,
so by $ \xCf (iM)$ $ \xbS' \xCK (\xbS'' - \xbG'') \xcc \xbS' \xCK
\xbS'' $ is small, so by $ \xCf (Opt)$ and $(s*s)$
$ \xbS'' - \xbG'' \xcc \xbS'' $ is small, so $ \xbG'' \xcc \xbS'' $
is big, and likewise $ \xbG' \xcc \xbS' $
is big.

$ \xcz $
\\[3ex]

\bd

$\hspace{0.01em}$


\label{Definition Mul-GH-Rel}

Call a relation $ \xeb $ a $ \xCf GH$ $(=$ general Hamming) relation iff
the following two
conditions hold:

$ \xCf (GH1)$ $ \xbs \xec \xbt $ $ \xcu $ $ \xbs' \xec \xbt' $ $ \xcu $
$(\xbs \xeb \xbt $ $ \xco $ $ \xbs' \xeb \xbt')$ $ \xch $ $ \xbs \xbs
' \xeb \xbt \xbt' $

(where $ \xbs \xec \xbt $ iff $ \xbs \xeb \xbt $ or $ \xbs = \xbt)$

$ \xCf (GH2)$ $ \xbs \xbs' \xeb \xbt \xbt' $ $ \xch $ $ \xbs \xeb \xbt $
$ \xco $ $ \xbs' \xeb \xbt' $

$ \xCf (GH2)$ means that some compensation is possible, e.g., $ \xbt \xeb
\xbs $ might be the
case, but $ \xbs' \xeb \xbt' $ wins in the end, so $ \xbs \xbs' \xeb
\xbt \xbt'.$

We use $ \xCf (GH)$ for $(GH1)+(GH2).$

\ed

\be

$\hspace{0.01em}$


\label{Example Mul-GH-Rel}

The following are examples of $ \xCf GH$ relations:

Define on all components $X_{i}$ a relation $ \xeb_{i}.$

(1) The set variant Hamming relation:

Let the relation $ \xeb $ be defined on $ \xbP \{X_{i}:i \xbe I\}$ by $
\xbs \xeb \xbt $ iff for all $j$ $ \xbs_{j} \xec_{j} \xbt_{j},$
and there is at least one $i$ s.t. $ \xbs_{i} \xeb_{i} \xbt_{i}.$

(2) The counting variant Hamming relation:

Let the relation $ \xeb $ be defined on $ \xbP \{X_{i}:i \xbe I\}$ by $
\xbs \xeb \xbt $ iff the number of $i$
such that $ \xbs_{i} \xeb_{i} \xbt_{i}$ is bigger than the number of $i$
such that $ \xbt_{i} \xeb_{i} \xbs_{i}.$

(3) The weighed counting Hamming relation:

Like the counting relation, but we give different (numerical) importance
to different $i.$ E.g., $ \xbs_{1} \xeb \xbt_{1}$ may count 1, $ \xbs_{2}
\xeb \xbt_{2}$ may count 2, etc.

$ \xcz $
\\[3ex]

\ee

\bp

$\hspace{0.01em}$


\label{Proposition Mul-GH-Rep}

Let $ \xbs \xeb \xbt \xcj \xbt \xce \xbm (\{ \xbs, \xbt \})$ and $ \xeb $
be smooth.
Then $ \xbm $ satisfies $(\xbm *1)$ (or, by
Fact \ref{Fact Mul-Big-Small} (page \pageref{Fact Mul-Big-Small})  equivalently
$(s*s))$
iff $ \xeb $ is a $ \xCf GH$ relation.

\ep

\subparagraph{
Proof
}

$\hspace{0.01em}$


(1) $(\xbm *1)$ entails the $ \xCf GH$ relation conditions

$ \xCf (GH1):$
Suppose $ \xbs \xeb \xbt $ and $ \xbs' \xec \xbt'.$ Then $ \xbt \xce
\xbm (\{ \xbs, \xbt \})=\{ \xbs \},$ and
$ \xbm (\{ \xbs', \xbt' \})=\{ \xbs' \}$ (either $ \xbs' \xeb \xbt'
$ or $ \xbs' = \xbt',$ so in both cases
$ \xbm (\{ \xbs', \xbt' \})=\{ \xbs' \}).$ As $ \xbt \xce \xbm (\{
\xbs, \xbt \}),$
$ \xbt \xbt' \xce \xbm (\{ \xbs, \xbt \} \xCK \{ \xbs', \xbt' \})$
$=_{(\xbm *1)}$ $ \xbm (\{ \xbs, \xbt \}) \xCK \xbm (\{ \xbs', \xbt'
\})$
$=\{ \xbs \} \xCK \{ \xbs' \}=\{ \xbs \xbs' \},$ so by smoothness $ \xbs
\xbs' \xeb \xbt \xbt'.$

$ \xCf (GH2):$
Let $X:=\{ \xbs, \xbt \},$ $Y:=\{ \xbs', \xbt' \},$ so $X \xCK Y=\{
\xbs \xbs', \xbs \xbt', \xbt \xbs', \xbt \xbt' \}.$
Suppose $ \xbs \xbs' \xeb \xbt \xbt',$ so $ \xbt \xbt' \xce \xbm (X
\xCK Y)$ $=_{(\xbm *1)}$ $ \xbm (X) \xCK \xbm (Y).$ If $ \xbs \xeB \xbt
,$ then
$ \xbt \xbe \xbm (X),$ likewise if $ \xbs' \xeB \xbt',$ then $ \xbt'
\xbe \xbm (Y),$ so $ \xbt \xbt' \xbe \xbm (X \xCK Y),$
contradiction.

(2) The $ \xCf GH$ relation conditions generate $(\xbm *1).$

$ \xbm (X \xCK Y) \xcc \xbm (X) \xCK \xbm (Y):$
Let $ \xbt \xbe X,$ $ \xbt' \xbe Y,$ $ \xbt \xbt' \xce \xbm (X) \xCK
\xbm (Y),$ then $ \xbt \xce \xbm (X)$ or $ \xbt' \xce \xbm (Y).$ Suppose
$ \xbt \xce \xbm (X),$ let $ \xbs \xbe X,$ $ \xbs \xeb \xbt,$ so by
condition $ \xCf (GH1)$ $ \xbs \xbt' \xeb \xbt \xbt',$ so
$ \xbt \xbt' \xce \xbm (X \xCK Y).$

$ \xbm (X) \xCK \xbm (Y) \xcc \xbm (X \xCK Y):$
Let $ \xbt \xbe X,$ $ \xbt' \xbe Y,$ $ \xbt \xbt' \xce \xbm (X \xCK Y),$
so there is $ \xbs \xbs' \xeb \xbt \xbt',$ $ \xbs \xbe X,$ $ \xbs'
\xbe Y,$
so by $ \xCf (GH2)$ either $ \xbs \xeb \xbt $ or $ \xbs' \xeb \xbt',$
so $ \xbt \xce \xbm (X)$ or $ \xbt' \xce \xbm (Y),$ so
$ \xbt \xbt' \xce \xbm (X) \xCK \xbm (Y).$

$ \xcz $
\\[3ex]

\bfa

$\hspace{0.01em}$


\label{Fact Mul-Small-Or}

(1) Let $ \xbG \xcc \xbS,$ $ \xbG' \xcc \xbS',$ $ \xbG \xCK \xbG'
\xcc \xbS \xCK \xbS' $ be small, let $ \xCf (GH2)$ hold, then
$ \xbG \xcc \xbS $ is small or $ \xbG' \xcc \xbS' $ is small.

(2) Let $ \xbG \xcc \xbS $ be small, $ \xbG' \xcc \xbS',$ let $ \xCf
(GH1)$ hold, then
$ \xbG \xCK \xbG' \xcc \xbS \xCK \xbS' $ is small.

\efa

\subparagraph{
Proof
}

$\hspace{0.01em}$


(1) Suppose $ \xbG \xcc \xbS $ is not small, so there is $ \xbg \xbe \xbG
$ and no $ \xbs \xbe \xbS $ with $ \xbs \xeb \xbg.$
Fix this $ \xbg.$
Consider $\{ \xbg \} \xCK \xbG'.$ As $ \xbG \xCK \xbG' \xcc \xbS \xCK
\xbS' $ is small, there is for each
$ \xbg \xbg',$ $ \xbg' \xbe \xbG' $ some $ \xbs \xbs' \xbe \xbS \xCK
\xbS',$ $ \xbs \xbs' \xeb \xbg \xbg'.$ By $ \xCf (GH2)$
$ \xbs \xeb \xbg $ or $ \xbs' \xeb \xbg',$
but $ \xbs \xeb \xbg $ was excluded, so for all $ \xbg' \xbe \xbG' $
there is $ \xbs' \xbe \xbS' $ with
$ \xbs' \xeb \xbg',$ so $ \xbG' \xcc \xbS' $ is small.

(2) Let $ \xbg \xbe \xbG,$ so there is $ \xbs \xbe \xbS $ and $ \xbs \xeb
\xbg.$ By $ \xCf (GH1),$ for any
$ \xbg' \xbe \xbG' $ $ \xbs \xbg' \xeb \xbg \xbg',$ so no $ \xbg \xbg
' \xbe \xbG \xCK \xbG' $ is minimal.

$ \xcz $
\\[3ex]

To complete our picture, we repeat from  \cite{GS09c} the
following very (perhaps too much so - see the discussion there) strong
definition and two results (the reader is referred there for proofs):

\bd

$\hspace{0.01em}$


\label{Definition Mul-Pr}

$(GH+)$ $ \xbs \xec \xbt $ $ \xcu $ $ \xbs' \xec \xbt' $ $ \xcu $ $(
\xbs \xeb \xbt $ $ \xco $ $ \xbs' \xeb \xbt')$ $ \xcj $ $ \xbs \xbs'
\xeb \xbt \xbt'.$

(Of course, $(GH+)$ entails $ \xCf (GH).)$

\ed

\bfa

$\hspace{0.01em}$


\label{Fact Mul-Pr}

$(\xbm *1)$ and $(\xbm *2)$ and the usual axioms for smooth relations
characterize
relations satisfying $(GH+).$

\efa

\bp

$\hspace{0.01em}$


\label{Proposition Mul-Pr}

Interpolation of the form $ \xbf \xcl \xba \xcn \xbq $ exists, if $(\xbm
*1)$ and $(\xbm *2)$ hold.

\ep

\paragraph{
Note
}

$\hspace{0.01em}$


\label{Section Note}

Note that already $(\xbm *1)$ results in a strong independence result in
the second
scenario:
Let $ \xbs \xbr' \xeb \xbt \xbr',$ then $ \xbs \xbr'' \xeb \xbt \xbr
'' $ for all $ \xbr''.$ Thus, whether $\{ \xbr'' \}$ is small,
or medium size (i.e. $ \xbr'' \xbe \xbm (\xbS')),$ the behaviour of $
\xbS \xCK \{ \xbr'' \}$ is the same.
This we do not have in the first scenario, as small sets may behave very
differently from medium size sets. (But, still, their internal structure
is the same, only the minimal elements change.)
When $(\xbm *2)$ holds, then if $ \xbs \xbs' \xeb \xbt \xbt' $ and $
\xbs \xEd \xbt,$ then $ \xbs \xeb \xbt,$ i.e. we
need not have $ \xbs' = \xbt'.$
\subsection{
Semantical interpolation
}

\label{Section Mul-Int}
\subsubsection{
Monotonic interpolation
}

\label{Section Mul-Mon-Int}

$ \xCO $

\label{Gin-Nota-Defin-3}
\begin{table}
\caption{Notation and Definitions}

\label{Table Gin-Not-Def}
\begin{center}
\tabcolsep=0.5pt.
\begin{tabular}{|c|c|c|}
\hline
\multicolumn{3}{|c|}{\bf Notation and definitions}\\
\hline

 \xEH 2-valued $\{0,1\}$ \xEH many-valued $(V, \xck)$ \xEP
\hline

language $L' \xcc L$ \xEH
\multicolumn{2}{|c|}{

propositional variables $s, \Xl $
} \xEP

semantic equivalence of $ \xbf,$ $ \xbq $ \xEH
\multicolumn{2}{|c|}{

$f_{ \xbf }=f_{ \xbq }$ (or for all $m$ $f_{m, \xbf }=f_{m, \xbq })$
} \xEP
\hline

definability of $f$ \xEH
\multicolumn{2}{|c|}{

$ \xcE \xbf:f_{ \xbf }=f$ (or for all $m$ $f_{m, \xbf }=f_{m})$
} \xEP
\hline

$ \xbG \xex L' $ \xEH
\multicolumn{2}{|c|}{

(for $ \xbG \xcc M)$ $ \xbG \xex L':=\{m \xex L':m \xbe \xbG \}$
} \xEP
\hline

model $m$ \xEH $m:L \xcp \{0,1\}$ \xEH $m:L \xcp V$ \xEP
$M$ set of all $L-$models \xEH \xEH \xEP
\hline

$m \xex L' $ \xEH
\multicolumn{2}{|c|}{

like $m,$ but restricted to $L' $
} \xEP
\hline

$m \xCq_{L' }m' $ \xEH
\multicolumn{2}{|c|}{

$m \xCq_{L' }m' $ iff $ \xcA s \xbe L'.m(s)=m' (s)$
} \xEP
\hline

model set of formula $ \xbf $ \xEH $M(\xbf) \xcc M,$
$f_{ \xbf }:M \xcp \{0,1\}$ \xEH $f_{ \xbf }:M \xcp V$ \xEP
\hline

general model set \xEH $M \xcc M,$
$f:M \xcp \{0,1\}$ \xEH $f:M \xcp V$ \xEP
\hline

$f$ insensitive to $L' $ \xEH
\multicolumn{2}{|c|}{

$ \xcA m,m' \xbe M.(m \xCq_{L-L' }m' \xch f(m)=f(m'))$
} \xEP
\hline

$f^{+}(m \xex L'),$ $f^{-}(m \xex L')$ \xEH
\multicolumn{2}{|c|}{

$f^{+}(m \xex L')=max\{f(m'):m' \xbe M,$ $m \xCq_{L' }m' \}$
} \xEP

 \xEH
\multicolumn{2}{|c|}{

$f^{-}(m \xex L')=min\{f(m'):m' \xbe M,$ $m \xCq_{L' }m' \}$
} \xEP
\hline

$f \xck g$ \xEH
\multicolumn{2}{|c|}{

$ \xcA m \xbe M.f(m) \xck g(m)$
} \xEP
\hline
\end{tabular}
\end{center}
\end{table}

$ \xCO $

$ \xCO $

\bd

$\hspace{0.01em}$


\label{Definition Gin-Base}

Let $M$ be the set of models for some language $ \xdl $ with set $L$ of
propositional
variables. Let $(V, \xck)$ be a finite, totally ordered set (of values).
Let $ \xbG \xcc M.$ $m,n$ etc. will be elements of $M.$ As usual, $ \xex $
will denote
the restriction of a function to part of its domain, and, by abuse of
language,
the restrictions of a set of functions.

(1) Let $J \xcc L,$ $f: \xbG \xcp V.$ Define
$f^{+}(m \xex J):=max\{f(m'):m \xex J=m' \xex J\}$ and
$f^{-}(m \xex J):=min\{f(m'):m \xex J=m' \xex J\}.$ (Similarly, if $m$ is
defined only
on $J,$ the condition is $m' \xex J=m,$ instead of $m \xex J=m' \xex J.)$

(2) Call $ \xbG $ rich iff for all $m,m' \xbe \xbG,$ $J \xcc L$ $(m \xex
J) \xcv (m' \xex (L-J)) \xbe \xbG.$
(I.e., we may cut and paste models.)

(3) Call $f: \xbG \xcp V$ insensitive to $J \xcc L$ iff for all $m,n$ $m
\xex (L-J)=n \xex (L- \xCf J)$
implies $f(m)=f(n)$ - i.e., the values of $m$ on $J$ have no importance
for $f.$

\ed

Let $L=J \xcv J' \xcv J'' $ be a disjoint union. If $f:M \xcp V$ is
insensitive to $J \xcv J'',$ we
can define for $m_{J' }:J' \xcp V$ $f(m_{J' })$ as any $f(m')$ such that
$m' \xex J' =m_{J' }.$

\bfa

$\hspace{0.01em}$


\label{Fact Gin-Pr-Int}

Let $ \xbG $ be rich, $f,g: \xbG \xcp V,$ $f(m) \xck g(m)$ for all $m \xbe
\xbG.$
Let $L=J \xcv J' \xcv J'',$ let $f$ be insensitive to $J,$ $g$ be
insensitive to $J''.$

Then $f^{+}(m_{J' }) \xck g^{-}(m_{J' })$ for all $m_{J' } \xbe \xbG \xex
J',$ and any $h: \xbG \xex J' \xcp V$ which is
insensitive to $J \xcv J'' $ is an interpolant iff

$f^{+}(m_{J' }) \xck h(m_{J' }) \xck g^{-}(m_{J' })$ for all $m_{J' } \xbe
\xbG \xex J'.$

$h$ can then be extended to the full $ \xbG $ in a unique way, as it is
insensitive to $J \xcv J''.$

\efa

\subparagraph{
Proof
}

$\hspace{0.01em}$


Let $L=J \xcv J' \xcv J'' $ be a pairwise disjoint union. Let $f$ be
insensitive to $J,$ $g$ be
insensitive to $J''.$

$h: \xbG \xcp V$ will have to be insensitive to $J \xcv J'',$ so we will
have to define
$h$ on $ \xbG \xex J',$ the extension to $ \xbG $ is then trivial.

Fix arbitrary $m_{J' }:J' \xcp V,$ $m_{J' }=m \xex J' $ for some $m \xbe
\xbG.$
We have $f^{+}(m_{J' }) \xck g^{-}(m_{J' }).$

Proof:
Choose $m_{J'' }$ such that $f^{+}(m_{J' })=f(m_{J}m_{J' }m_{J'' })$ for
any $m_{J}.$
(Recall that $f$ is insensitive to $J.)$
Let $n_{J'' }$ be one such $m_{J'' }.$
Likewise,
choose $m_{J}$ such that $g^{-}(m_{J' })=g(m_{J}m_{J' }m_{J'' })$ for any
$m_{J'' }.$
Let $n_{J}$ be one such $m_{J}.$
Consider $n_{J}m_{J' }n_{J'' } \xbe \xbG $ (recall that $ \xbG $ is rich).
By definition, $f^{+}(m_{J' })=f(n_{J}m_{J' }n_{J'' })$ and $g^{-}(m_{J'
})=g(n_{J}m_{J' }n_{J'' }),$ but
by prerequisite $f(n_{J}m_{J' }n_{J'' }) \xck g(n_{J}m_{J' }n_{J'' }),$ so
$f^{+}(m_{J' }) \xck g^{-}(m_{J' }).$

Thus, any $h$ such that $h$ is insensitive to $J \xcv J'' $ and

(Int) $f^{+}(m_{J' }) \xck h(m_{J' }) \xck g^{-}(m_{J' })$

is an interpolant for $f$ and $g.$

But (Int) is also a necessary condition.

Proof:

Suppose $h$ is insensitive to $J \xcv J'' $ and $h(m_{J' })<f^{+}(m_{J'
}).$ Let $n_{J'' }$ be as above, i.e.,
$f(m_{J}m_{J' }n_{J'' })=f^{+}(m_{J' })$ for any $m_{J}.$ Then
$h(m_{J}m_{J' }n_{J'' })=h(m_{J' })<f^{+}(m_{J' })=f(m_{J}m_{J' }n_{J''
}),$ so $h$ is not an interpolant.

The proof that $h(m_{J' })$ has to be $ \xck g^{-}(m_{J' })$ is analogous.

We summarize:

$f$ and $g$ have an interpolant $h,$ and
$h$ is an interpolant for $f$ and $g$ iff $h$ is insensitive to $J \xcv
J'' $ and
for any $m_{J' } \xbe \xbG \xex J' $
$f^{+}(m_{J' }) \xck h(m_{J' }) \xck g^{-}(m_{J' }).$

$ \xcz $
\\[3ex]

$ \xCO $
\subsubsection{
Non-monotonic interpolation
}

\label{Section Mul-NonMon-Int}

\bp

$\hspace{0.01em}$


\label{Proposition Mul-Mu*1-Int}

$(\xbm *1)$ entails semantical interpolation of the form $ \xbf \xcn \xba
\xcn \xbq $ in 2-valued
non-monotonic logic generated by minimal model sets.
(As the model sets might not be definable, syntactic interpolation does
not follow automatically.)

\ep

\subparagraph{
Proof
}

$\hspace{0.01em}$


Let the product be defined on $J \xcv J' \xcv J'' $ (i.e., $J \xcv J' \xcv
J'' $ is the set
of propositional variables in the intended application).
Let $ \xbf $ be defined on $J' \xcv J'',$ $ \xbq $ on $J \xcv J'.$
See Diagram \ref{Diagram Mul-Base-2} (page \pageref{Diagram Mul-Base-2}).

We abuse notation and write $ \xbf \xcn \xbS $ if $ \xbm (\xbf) \xcc
\xbS.$ As usual,
$ \xbm (\xbf)$ abbreviates $ \xbm (M(\xbf)).$

For clarity, even if it clutters up notation, we will be precise about
where
$ \xbm $ is formed. Thus, we write $ \xbm_{J \xcv J' \xcv J'' }(X)$ when
we take the minimal elements
in the full product, $ \xbm_{J}(X)$ when we consider only the product on
$J,$ etc.

Let $ \xbf \xcn \xbq,$ i.e., $ \xbm_{J \xcv J' \xcv J'' }(\xbf) \xcc M(
\xbq).$
We show that $X_{J} \xCK (\xbm_{J \xcv J' \xcv J'' }(\xbf) \xex J')
\xCK X_{J'' },$ i.e., that
$ \xbm_{J \xcv J' \xcv J'' }(\xbf)$ $ \xcc $ $X_{J} \xCK (\xbm_{J \xcv
J' \xcv J'' }(\xbf) \xex J') \xCK X_{J'' },$
and that
$ \xbm_{J \xcv J' \xcv J'' }(X_{J} \xCK (\xbm_{J \xcv J' \xcv J'' }(\xbf
) \xex J') \xCK X_{J'' }) \xcc M(\xbq).$

The first property is trivial, we turn to the second.
(1) As $M(\xbf)=X_{J} \xCK M(\xbf) \xex (J' \xcv J''),$
$ \xbm_{J \xcv J' \xcv J'' }(\xbf)$ $=$ $ \xbm_{J}(X_{J}) \xCK \xbm_{J'
\xcv J'' }(M(\xbf) \xex (J' \xcv J''))$ by $(\xbm *1).$

(2) By $(\xbm *1),$
$ \xbm_{J \xcv J' \xcv J'' }(X_{J} \xCK (\xbm_{J \xcv J' \xcv J'' }(\xbf
) \xex J') \xCK X_{J'' })$ $=$
$ \xbm_{J}(X_{J}) \xCK \xbm_{J' }(\xbm_{J \xcv J' \xcv J'' }(\xbf) \xex
J') \xCK \xbm_{J'' }(X_{J'' }).$

So it suffices to show $ \xbm_{J}(X_{J}) \xCK \xbm_{J' }(\xbm_{J \xcv J'
\xcv J'' }(\xbf) \xex J') \xCK \xbm_{J'' }(X_{J'' }) \xcm \xbq.$

Proof:
Let $ \xbs = \xbs_{J} \xbs_{J' } \xbs_{J'' }$ $ \xbe $
$ \xbm_{J}(X_{J}) \xCK \xbm_{J' }(\xbm_{J \xcv J' \xcv J'' }(\xbf) \xex
J') \xCK \xbm_{J'' }(X_{J'' }),$
so $ \xbs_{J} \xbe \xbm_{J}(X_{J}).$

By $ \xbm_{J' }(\xbm_{J \xcv J' \xcv J'' }(\xbf) \xex J') \xcc \xbm_{J
\xcv J' \xcv J'' }(\xbf) \xex J',$
there is $ \xbs' = \xbs_{J}' \xbs_{J' }' \xbs_{J'' }' \xbe \xbm_{J \xcv
J' \xcv J'' }(\xbf)$
s.t. $ \xbs'_{J' }= \xbs_{J' },$ i.e. $ \xbs' = \xbs'_{J} \xbs_{J' }
\xbs'_{J'' }.$ As $ \xbs' \xbe \xbm_{J \xcv J' \xcv J'' }(\xbf),$ $
\xbs' \xcm \xbq.$

By (1) and $ \xbs_{J} \xbe \xbm_{J}(X_{J})$ also $ \xbs_{J} \xbs_{J' }
\xbs'_{J'' } \xbe \xbm_{J \xcv J' \xcv J'' }(\xbf),$ so also
$ \xbs_{J} \xbs_{J' } \xbs'_{J'' } \xcm \xbq.$

But $ \xbq $ does not depend on $J'',$ so also
$ \xbs $ $=$ $ \xbs_{J} \xbs_{J' } \xbs_{J'' }$ $ \xcm $ $ \xbq.$

$ \xcz $
\\[3ex]

\br

$\hspace{0.01em}$


\label{Remark Mul-Many}

We can try to extend our result to many-vanued logics. But the we have
first
to make precise what we want to do. One approach might be:
Let max be the maximal truth value. We look at the set of models where
a formula has truth value max, and then look at the minimal models of this
set, under some relation. But we can also consider other ideas: we can
look
at all truth values separately, do minimization for all values separately,
etc.

\er

$ \xCO $

\vspace{10mm}

\begin{diagram}

\label{Diagram Mul-Base-2}
\index{Diagram Mul-Base-2}

\centering
\setlength{\unitlength}{1mm}
{\renewcommand{\dashlinestretch}{30}
\begin{picture}(150,100)(0,0)

\path(20,90)(90,90)(90,30)(20,30)(20,90)

\path(20,89.5)(40,89.5)
\path(70,89.5)(90,89.5)
\path(40,49.5)(70,49.5)

\path(20,30)(20,28)
\path(90,30)(90,28)

\path(40,90)(40,28)
\path(70,90)(70,28)

\path(20,70)(70,70)
\path(40,80)(90,80)
\path(20,50)(90,50)

\put(93,80){{\xssc $\xbf$}}
\put(16,70){{\xssc $\xbq$}}
\put(93,50){{\xssc $\xbm(\xbf)$}}

\put(30,24){{\xssc $J$}}
\put(54,24){{\xssc $J'$}}
\put(80,24){{\xssc $J"$}}

\put(30,14){{\xssc Non-monotonic interpolation}}
\put(31,10){{\xssc Double lines: interpolant}}

\end{picture}
}

\end{diagram}

\vspace{4mm}

$ \xCO $
\paragraph{
Remarks for the converse: from interpolation to $(\xbm *1)$
}

\be

$\hspace{0.01em}$


\label{Example Mul-Mu}

We show here in (1) and (2) that half of the condition $(\xbm *1)$ is not
sufficient for interpolation, and in (3) that interpolation may hold, even
if
$(\xbm *1)$ fails. When looking closer, the latter is not surprising: $
\xbm $ of
sub-products may be defined in a funny way, which has nothing to do with
the
way $ \xbm $ on the big product is defined.

Consider the language based on $p,q,r.$

For (1) and (2) define the order $ \xeb $ on sequences of length 3 by
$ \xCN p \xCN q \xCN r \xeb p \xCN q \xCN r,$ leave all other 3-sequences
incomparabel.

Let $ \xbf = \xCN q \xcu \xCN r,$ $ \xbq = \xCN p \xcu \xCN q,$ so $ \xbm
(\xbf)= \xCN p \xcu \xCN q \xcu \xCN r,$ and $ \xbf \xcn \xbq.$
Suppose there is $ \xba,$ $ \xbf \xcn \xba \xcn \xbq,$ $ \xba $ written
with $q$ only, so $ \xba $ is equivalent
to FALSE, TRUE, $q,$ or $ \xCN q.$ $ \xbf \xcN FALSE,$ $ \xbf \xcN q.$
$TRUE \xcN \xbq,$ $ \xCN q \xcN \xbq.$ Thus, there
is no such $ \xba,$ and $ \xcn $ has no interpolation. We show in (1) and
(2) that
we can make both directions of $(\xbm *1)$ true separately, so they do
not
suffice to obtain interpolation.

(1) We make $ \xbm (X \xCK Y) \xcc \xbm (X) \xCK \xbm (Y)$ true, but not
the converse.

Do not order any sequences of length 2 or 1, i.e. $ \xbm $ is there always
identity. Thus, $ \xbm (X \xCK Y) \xcc X \xCK Y= \xbm (X) \xCK \xbm (Y)$
holds trivially.

For (2) and (3), consider the following ordering $<$ between sequences:
$ \xbs < \xbt $ iff there is $ \xCN x$ is $ \xbs,$ $x$ in $ \xbt,$ but
for no $y$ $y$ in $ \xbs,$ $ \xCN y$ in $ \xbt.$
E.g., $ \xCN p<p,$ $ \xCN pq<pq,$ $ \xCN p \xCN q<pq,$ but $ \xCN pq \xEc
p \xCN q.$

(2) We make $ \xbm (X \xCK Y) \xcd \xbm (X) \xCK \xbm (Y)$ true, but not
the converse.

We order all sequences of length 1 or 2 by $<.$

Suppose $ \xbs \xbe X \xCK Y- \xbm (X \xCK Y).$ Case 1: $X \xCK Y$
consists of sequences of length 2.
Then, by definition, $ \xbs \xce \xbm (X) \xCK \xbm (Y).$ Case 2: $X \xCK
Y$ consists of sequences
of length 3. Then $ \xbs =p \xCN q \xCN r,$ and there is $ \xbt = \xCN p
\xCN q \xCN r \xbe X \xCK Y.$ So
$\{p, \xCN p\} \xcc X$ or $\{p \xCN q, \xCN p \xCN q\} \xcc X,$ but in
both cases $ \xbs \xex X \xce \xbm (X).$

Finally, note that $ \xbm (TRUE) \xcC \{ \xCN p \xCN q \xCN r\},$ so full
$(\xbm *1)$ does not hold.

(3) We make interpolation hold, but $ \xbm (X) \xCK \xbm (Y) \xcC \xbm (X
\xCK Y):$

We order all sequences of length 3 by $<.$ Shorter sequences are made
incomparabel, so for shorter sequences $ \xbm (X)=X.$

Obviously, in general $ \xbm (X) \xCK \xbm (Y) \xcC \xbm (X \xCK Y).$

But the proof of
Proposition \ref{Proposition Mul-Mu*1-Int} (page \pageref{Proposition
Mul-Mu*1-Int})  goes through as above, only
directly, without the
use of factorizing and taking $ \xbm $ of the factors.

$ \xcz $
\\[3ex]
\subsection{
Language change
}

\ee

Independence of language fragments gives us the following perspectives:

 \xEh

 \xDH it makes independent and parallel treatment of fragments possible,
and offers thus efficient treatment in applications
(descriptive logics etc.).

 \xDH it results in new rules similar to the classical ones like AND, OR,
Cumulativity, etc. We can thus obtain postulates about reasonable
behaviour,
but also classification by those rules, see
Table \ref{Table Mul-Laws} (page \pageref{Table Mul-Laws}), Scenario 2, Logical
property.

 \xDH it sheds light on notions like
``ceteris paribus'', which we saw in the context of obligations,
see  \cite{GS08g}.

 \xDH it clarifies notions like
``normal with respect to $ \xbf,$ but not $ \xbq $''

 \xDH it helps to understand e.g. inheritance diagrams where
arrows make other information accessible, and we need an underlying
mechanism to combine bits of information, given in different
languages.

 \xEj
\subsection{
A relevance problem
}

\label{Section Mul-Relev}

Consider the formula $ \xbf:=a \xcu \xCN a \xcu b.$ Then $M(\xbf)= \xCQ
.$ But we cannot recover
where the problem came from, and this results in the EFQ rule. We now
discuss
one, purely algebraic, approach to remedy.

Consider 3 valued models, with a new value $b$ for both, in addition to
$t$ and $f.$
Above formula would then have the model $m(a)=b,$ $m(b)=t.$ So there is a
model,
EFQ fails, and we can recover the culprit.

To have the usual behaviour of $ \xcu $ as intersection, it might be good
to change
the definition so that $m(x)=b$ is always a model. Then $M(b)=\{m(b)=t,m'
(b)=b\},$
$M(\xCN b)=\{m(b)=f,m' (b)=b\},$ and $M(b \xcu \xCN b)=\{m' (b)=b\}.$

It is not yet clear which version to choose, and we have no syntactic
characterization.

Other idea:

Use meaningless models. Take a conjunction of literals.
$m(a)=t$ and $m(a)=x$ is a model if there is only a in the conjunction,
$m(a)=f$ and $m(a)=x$ if there is only $ \xCN a$ in the conjunction,
$m(a)=*$ if both are present, and all models, if none is present.
Thus there is always a model, and we can isolate the contradictory parts:
there, only $m(a)=x$ is present.
\subsection{
Small subspaces
}

When considering small subsets in nonmonotonic logic, we neglect small
subsets of models.
What is the analogue when considering small subspaces, i.e. when
$J=J' \xcv J'',$ with $J'' $ small in $J$ in nonmonotonic logic?

It is perhaps easiest to consider the relation based approach first.
So we have an order on $ \xbP J' $ and one on $ \xbP J'',$ $J'' $ is
small, and we want
to know how to construct a corresponding order on $ \xbP J.$ Two solutions
come
to mind:

 \xEI

 \xDH a less radical one: we make a lexicographic ordering, where the one
on $ \xbP J' $
has precedence over the one on $ \xbP J'',$

 \xDH a more radical one: we totally forget about the ordering of $ \xbP
J'',$ i.e. we
do as if the ordering on $ \xbP J'' $ were the empty set, i.e.
$ \xbs' \xbs'' \xeb \xbt' \xbt'' $ iff $ \xbs' \xeb \xbt' $ and $
\xbs'' = \xbt''.$

We call this condition $forget(J'').$

 \xEJ

The less radical one is already covered by our relation conditions
$ \xCf (GH).$
The more radical one is probably more interesting. Suppose $ \xbf' $ is
written in
language $J',$ $ \xbf'' $ in language $J'',$ we then have

$ \xbf' \xcu \xbf'' \xcn \xbq' \xcu \xbq'' $ iff $ \xbf' \xcn \xbq'
$ and $ \xbf'' \xcl \xbq''.$

This approach, is of course the same as considering on the small
coordinate only ALL as a big subset, (see the lines $x*1/1*x$
in Table \ref{Table Mul-Laws} (page \pageref{Table Mul-Laws})),
so, in principle, we get nothing new.
\section{
Revision and distance relations
}

We will look here into distance based theory revision a la
AGM, see  \cite{AGM85} and  \cite{LMS01}, and
also  \cite{Sch04} for more details.
First, we introduce some notation, and give a result taken from
 \cite{GS09c} (slightly modified).

\bd

$\hspace{0.01em}$


\label{Definition Mul-GHD}

Let $d$ be a distance on some product space $X \xCK Y,$ and its
components.
(We require of distances only that they are comparable, that $d(x,y)=0$
iff
$x=y,$ and that $d(x,y) \xcg 0.)$

$d$ is called a generalized Hamming distance $ \xCf (GHD)$
iff it satisfies the following two properties:

$ \xCf (GHD1)$ $d(\xbs, \xbt) \xck d(\xba, \xbb)$ and $d(\xbs',
\xbt') \xck d(\xba', \xbb')$ and
$(d(\xbs, \xbt)<d(\xba, \xbb)$ or $d(\xbs', \xbt')<d(\xba',
\xbb'))$ $ \xch $ $d(\xbs \xbs', \xbt \xbt')<d(\xba \xba', \xbb
\xbb')$

$ \xCf (GHD2)$ $d(\xbs \xbs', \xbt \xbt')<d(\xba \xba', \xbb \xbb
')$ $ \xch $
$d(\xbs, \xbt)<d(\xba, \xbb)$ or $d(\xbs', \xbt')<d(\xba',
\xbb')$

\ed

\bd

$\hspace{0.01em}$


\label{Definition Sin-Bar}

Given a distance $d,$ define for two sets $X,Y$

$X \xfA Y:=\{y \xbe Y: \xcE x \xbe X(\xCN \xcE x' \xbe X,y' \xbe Y.d(x'
,y')<d(x,y))\}.$

We assume that $X \xfA Y \xEd \xCQ $ if $X,Y \xEd \xCQ.$ Note that this
is related to the
consistency axiom of AGM theory revision: revising by a consistent formula
gives a consistent result. The assumption may be wrong due to infinite
descending chains of distances.

\ed

\bd

$\hspace{0.01em}$


\label{Definition Sin-TR}

Given $ \xfA $ on models, we can define an AGM revision operator $*$ as
follows:

$T* \xbf:=Th(M(T) \xfA M(\xbf))$

where $T$ is a theory, and $Th(X)$ is the set of formulas which hold in
all $x \xbe X.$

It was shown in  \cite{LMS01} that a revision operator thus defined
satisfies the AGM revision postulates.

\ed

We have a result analogous to the relation case:

\bfa

$\hspace{0.01em}$


\label{Fact Mul-GHD-Bar}

Let $ \xfA $ be defined by a generalized Hamming distance, then $ \xfA $
satisfies

$(\xfA *)$ $(\xbS_{1} \xCK \xbS_{1}') \xfA (\xbS_{2} \xCK \xbS_{2}'
)=(\xbS_{1} \xfA \xbS_{2}) \xCK (\xbS_{1}' \xfA \xbS_{2}').$

\efa

\subparagraph{
Proof
}

$\hspace{0.01em}$


``$ \xcc $'':

Suppose $d(\xbs \xbs', \xbt \xbt')$ is minimal. If there is
$ \xba \xbe \xbS_{1},$ $ \xbb \xbe \xbS_{2}$ s.t. $d(\xba, \xbb)<d(
\xbs, \xbt),$ then
$d(\xba \xbs', \xbb \xbt')<d(\xbs \xbs', \xbt \xbt')$ by $ \xCf
(GHD1),$ so $d(\xbs, \xbt)$ and $d(\xbs', \xbt')$ have to be
minimal.

``$ \xcd $'':

For the converse, suppose $d(\xbs, \xbt)$ and $d(\xbs', \xbt')$
are minimal, but
$d(\xbs \xbs', \xbt \xbt')$ is not, so $d(\xba \xba', \xbb \xbb'
)<d(\xbs \xbs', \xbt \xbt')$ for some $ \xba \xba', \xbb \xbb',$
then $d(\xba, \xbb)<d(\xbs, \xbt)$ or $d(\xba', \xbb')<d(\xbs
', \xbt')$ by $ \xCf (GHD2),$ contradiction.

$ \xcz $
\\[3ex]

These properties translate to logic as follows:

\bco

$\hspace{0.01em}$


\label{Corollary Mul-GHD-TR}

If $ \xbf $ and $ \xbq $ are defined on a separate language from that of $
\xbf' $ and $ \xbq',$
and the distance satisfies $ \xCf (GHD1)$ and $ \xCf (GHD2),$ then for
revision holds:

$(\xbf \xcu \xbf')*(\xbq \xcu \xbq')=(\xbf * \xbq) \xcu (\xbf' *
\xbq').$
\subsection{
Interpolation for distance based revision
}

\eco

The limiting condition (consistency) imposes a strong restriction:
Even for $ \xbf *TRUE,$ the result may need many variables (those in $
\xbf).$

\bl

$\hspace{0.01em}$


\label{Lemma Mul-TR}

Let $ \xfA $ satisfy $(\xfA *).$

Let $J \xcc L,$ $ \xbr $ be written in sublanguage $J,$ let $ \xbf, \xbq
$ be written in $L- \xCf J,$
let $ \xbf', \xbq' $ be written in $J' \xcc J.$

Let $(\xbf \xcu \xbf')*(\xbq \xcu \xbq') \xcl \xbr,$ then $ \xbf'
* \xbq' \xcl \xbr.$

(This is suitable interpolation, but we also need to factorize the
revision construction.)

\el

\subparagraph{
Proof
}

$\hspace{0.01em}$


$(\xbf \xcu \xbf')*(\xbq \xcu \xbq')$ $=$ $(\xbf * \xbq) \xcu (
\xbf' * \xbq')$ by $(\xfA *).$
So $(\xbf \xcu \xbf')*(\xbq \xcu \xbq') \xcl \xbr $ iff $(\xbf *
\xbq) \xcu (\xbf' * \xbq') \xcl \xbr,$ but
$(\xbf * \xbq) \xcu (\xbf' * \xbq') \xcl \xbr $ iff $(\xbf' * \xbq
') \xcl \xbr,$ as $ \xbr $ contains no variables
of $ \xbf' $ or $ \xbq'.$
$ \xcz $
\\[3ex]
\section{
Summary of properties
}

\label{Section Mul-Sum}

$pr(b)=b$ means: the projection of a big set on one of its coordinates
is big again.
\begin{table}
\caption{Multiplication laws}

\label{Table Mul-Laws}
\begin{center}
\tabcolsep=0.5pt.
\begin{turn}{90}
{\xssB
\begin{tabular}{|c|c|c|c|c|c|c|c|c|}
\hline
\multicolumn{9}{|c|}{

Multiplication laws
} \xEP
\hline
Multiplication \xEH
\multicolumn{3}{|c|}{

Scenario 1
} \xEH
\multicolumn{5}{|c|}{

Scenario 2 $(*$ symmetrical, only 1 side shown)
} \xEP

law \xEH
\multicolumn{3}{|c|}{

(see Diagram \ref{Diagram Mul-Add} (page \pageref{Diagram Mul-Add}))
} \xEH
\multicolumn{5}{|c|}{

(see Diagram \ref{Diagram Mul-Prod} (page \pageref{Diagram Mul-Prod}))
} \xEP
\cline{2-9}

 \xEH Corresponding algebraic \xEH Logical property \xEH Relation \xEH
Algebraic property \xEH
Logical property \xEH
\multicolumn{3}{|c|}{

Interpolation
} \xEP
\cline{7-9}

 \xEH addition property \xEH \xEH property \xEH $(\xbG_{i} \xcc
\xbS_{i})$ \xEH
$ \xba, \xbb $ in $ \xdl_{1},$ $ \xba', \xbb' $ in $ \xdl_{2}$ \xEH
Multiplic. \xEH Relation \xEH Inter- \xEP
 \xEH \xEH \xEH \xEH \xEH
$ \xdl = \xdl_{1} \xcv \xdl_{2}$ (disjoint) \xEH law \xEH property \xEH
polation \xEP
\hline
\multicolumn{9}{|c|}{

Non-monotonic logic
} \xEP
\hline

$x*1 \xch x$ \xEH
\multicolumn{3}{|c|}{trivial}

 \xEH
$ \xbG_{1} \xbe \xdf (\xbS_{1})$ $ \xch $
 \xEH $ \xba \xcn_{ \xdl_{i}} \xbb $ $ \xch $ $ \xba \xcn_{ \xdl } \xbb $
\xEH \xEH \xEH \xEP
\cline{1-4}

$1*x \xch x$ \xEH
\multicolumn{3}{|c|}{trivial}

 \xEH
$ \xbG_{1} \xCK \xbS_{2} \xbe \xdf (\xbS_{1} \xCK \xbS_{2})$
 \xEH $ \xba \xcn_{ \xdl } \xbb $ \xEH \xEH \xEH \xEP
\hline

$x*s \xch s$ \xEH (iM) \xEH $ \xba \xcn \xCN \xbb $ $ \xch $ \xEH -
 \xEH
dual to $x*1 \xch 1$
 \xEH $ \xba \xcn_{ \xdl_{1}} \xbb $ $ \xch $ $ \xba \xcn_{ \xdl } \xbb $
\xEH \xEH \xEH \xEP
 \xEH
$A \xcc B \xbe \xdi (X)$ $ \xch $ $A \xbe \xdi (X)$
 \xEH $ \xba \xcn \xCN \xbb \xco \xbg $ \xEH \xEH \xEH \xEH \xEH \xEH \xEP
\cline{1-4}

$s*x \xch s$ \xEH $(eM \xdi)$ \xEH $ \xba \xcu \xbb \xcn \xCN \xbg $ $
\xch $ \xEH - \xEH \xEH \xEH \xEH \xEH \xEP
 \xEH
$X \xcc Y$ $ \xch $ $ \xdi (X) \xcc \xdi (Y),$
 \xEH $ \xba \xcn \xCN \xbb \xco \xCN \xbg $ \xEH \xEH \xEH \xEH \xEH \xEH
\xEP
 \xEH
$X \xcc Y$ $ \xch $
 \xEH \xEH \xEH \xEH \xEH \xEH \xEH \xEP
 \xEH
$ \xdf (Y) \xcs \xdp (X) \xcc \xdf (X)$
 \xEH \xEH \xEH \xEH \xEH \xEH \xEH \xEP
\hline

$b*b \xch b$ \xEH $(< \xbo *s),$ $(\xdm^{+}_{ \xbo })$ (3) \xEH $ \xba
\xcn \xbb,$ $ \xba \xcu \xbb \xcn \xbg $ \xEH - (Filter)
 \xEH
$ \xbG_{1} \xbe \xdf (\xbS_{1}), \xbG_{2} \xbe \xdf (\xbS_{2})$ $ \xch $
 \xEH $ \xba \xcn_{ \xdl_{1}} \xbb,$ $ \xba' \xcn_{ \xdl_{2}} \xbb' $ $
\xch $ \xEH $b*b \xcj b:$
 \xEH $ \xCf (GH)$ \xEH $ \xcn \xDO \xcn $ \xEP
 \xEH
$A \xbe \xdf (X),X \xbe \xdf (Y)$ $ \xch $
 \xEH $ \xch $ $ \xba \xcn \xbg $ \xEH \xEH
$ \xbG_{1} \xCK \xbG_{2} \xbe \xdf (\xbS_{1} \xCK \xbS_{2})$
 \xEH
$ \xba \xcu \xba' \xcn_{ \xdl } \xbb \xcu \xbb' $
 \xEH $(\xbm *1)$ \xEH \xEH \xEP
 \xEH $A \xbe \xdf (Y)$ \xEH \xEH \xEH \xEH \xEH $ \xcj $ \xEH \xEH \xEP
\cline{1-6}

$b*m \xch m$ \xEH $(< \xbo *s),$ $(\xdm^{+}_{ \xbo })$ (2) \xEH $ \xba
\xcN \xCN \xbb,$ $ \xba \xcu \xbb \xcn \xbg $
 \xEH - (Filter) \xEH
$ \xbG_{1} \xbe \xdf (\xbS_{1}), \xbG_{2} \xbe \xdm^{+}(\xbS_{2})$ $
\xch $
 \xEH
$ \xba \xcN_{ \xdl_{1}} \xCN \xbb,$ $ \xba' \xcn_{ \xdl_{2}} \xbb' $ $
\xch $
 \xEH $(s*s)$ \xEH \xEH \xEP
 \xEH
$A \xbe \xdm^{+}(X),X \xbe \xdf (Y)$ $ \xch $
 \xEH $ \xch $ $ \xba \xcN \xCN \xbb \xco \xCN \xbg $
 \xEH \xEH
$ \xbG_{1} \xCK \xbG_{2} \xbe \xdm^{+}(\xbS_{1} \xCK \xbS_{2})$
 \xEH
$ \xba \xcu \xba' \xcN_{ \xdl } \xCN \xbb \xco \xbb' $
 \xEH \xEH \xEH \xEP
 \xEH $A \xbe \xdm^{+}(Y)$ \xEH \xEH \xEH \xEH \xEH \xEH \xEH \xEP
\cline{1-4}

$m*b \xch m$ \xEH $(< \xbo *s),$ $(\xdm^{+}_{ \xbo })$ (1) \xEH $ \xba
\xcn \xbb,$ $ \xba \xcu \xbb \xcN \xCN \xbg $
 \xEH - (Filter) \xEH
 \xEH \xEH \xEH \xEH \xEP
 \xEH
$A \xbe \xdf (X),X \xbe \xdm^{+}(Y)$ $ \xch $
 \xEH $ \xch $ $ \xba \xcN \xCN \xbb \xco \xCN \xbg $
 \xEH \xEH \xEH \xEH \xEH \xEH \xEP
 \xEH $A \xbe \xdm^{+}(Y)$ \xEH \xEH \xEH \xEH \xEH \xEH \xEH \xEP
\cline{1-6}

$m*m \xch m$ \xEH $(\xdm^{++})$ \xEH Rational Monotony \xEH ranked \xEH
$ \xbG_{1} \xbe \xdm^{+}(\xbS_{1}), \xbG_{2} \xbe \xdm^{+}(\xbS_{2})$
 \xEH
$ \xba \xcN_{ \xdl_{1}} \xCN \xbb,$ $ \xba' \xcn_{ \xdl_{2}} \xCN \xbb'
$ $ \xch $
 \xEH \xEH \xEH \xEP
 \xEH
$A \xbe \xdm^{+}(X),X \xbe \xdm^{+}(Y)$
 \xEH \xEH \xEH $ \xch $
 \xEH
$ \xba \xcu \xba' \xcN_{ \xdl } \xCN \xbb \xco \xCN \xbb' $
 \xEH \xEH \xEH \xEP
 \xEH $ \xch $ $A \xbe \xdm^{+}(Y)$ \xEH \xEH \xEH
$ \xbG_{1} \xCK \xbG_{2} \xbe \xdm^{+}(\xbS_{1} \xCK \xbS_{2})$
 \xEH \xEH \xEH \xEH \xEP
\hline

$b*b \xcj b$ $+$ \xEH \xEH \xEH \xEH \xEH
$ \xba \xcn \xbb $ $ \xch $ $ \xba \xex \xdl_{1} \xcn \xbb \xex \xdl_{1}$
 \xEH \xEH $(GH+)$ \xEH $ \xcl \xDO \xcn $ \xEP
$pr(b)=b$ \xEH \xEH \xEH \xEH \xEH
and
 \xEH \xEH \xEH \xEP
 \xEH \xEH \xEH \xEH \xEH
$ \xba \xcn_{ \xdl_{1}} \xbb,$ $ \xba' \xcn_{ \xdl_{2}} \xbb' $ $ \xch
$
 \xEH \xEH \xEH \xEP
 \xEH \xEH \xEH \xEH \xEH
$ \xba \xcu \xba' \xcn_{ \xdl } \xbb \xcu \xbb' $
 \xEH \xEH \xEH \xEP
\hline

$J' $ small \xEH \xEH \xEH \xEH
 \xEH $ \xba \xcu \xba' \xcn \xbb \xcu \xbb' $ $ \xcj $
 \xEH \xEH $forget(J')$ \xEH - \xEP
 \xEH \xEH \xEH \xEH \xEH
$ \xba \xcn \xbb,$ $ \xba' \xcl \xbb' $ \xEH \xEH \xEH \xEP
\hline
\multicolumn{9}{|c|}{

Theory revision
} \xEP
\hline

 \xEH \xEH \xEH \xEH $(\xfA *):$ \xEH \xEH \xEH $ \xCf (GHD)$ \xEH $(
\xbf \xcu \xbf')*(\xbq \xcu \xbq') \xcl \xbr $ \xEP
 \xEH \xEH \xEH \xEH $(\xbS_{1} \xCK \xbS_{1}') \xfA (\xbS_{2} \xCK
\xbS_{2}')=$ \xEH $(\xbf \xcu \xbf')*(\xbq \xcu \xbq')=$
 \xEH \xEH \xEH $ \xch $ $ \xbf' * \xbq' \xcl \xbr $ \xEP
 \xEH \xEH \xEH \xEH $(\xbS_{1} \xfA \xbS_{2}) \xCK (\xbS_{1}' \xfA
\xbS_{2}')$ \xEH
$(\xbf * \xbq) \xcu (\xbf' * \xbq')$
 \xEH \xEH \xEH $ \xbf, \xbq $ in $J,$ \xEP
 \xEH \xEH \xEH \xEH
 \xEH \xEH \xEH \xEH $ \xbf', \xbq', \xbr $ in $L- \xCf J$ \xEP
\hline
\end{tabular}
}
\end{turn}
\end{center}
\end{table}

Note that $A \xCK B \xcc X \xCK Y$ big $ \xch $ $A \xcc X$ big etc. is
intuitively better justified
than the other direction, as the proportion might increase in the latter,
decrease in the former. Cf.
the table ``Rules on size'',
Section \ref{Section Mul-Intro} (page \pageref{Section Mul-Intro}),
``increasing proportions''.
\clearpage

\vspace{3mm}


\end{document}